\documentclass{amsart}
\usepackage{graphicx}
\usepackage{amsmath}
\usepackage{amssymb}
\usepackage{color}
\newcommand{\eps}{\varepsilon}
\newcommand{\p}{\partial}

\newcommand{\Ome}{\Omega}
\newcommand{\nab}{\nabla}
\newcommand{\G}{\Gamma}
\newcommand{\Del}{\Delta}

\newcommand{\bu}{\mathbf{u}}

\newcommand{\bv}{\mathbf{v}}
\newcommand{\bR}{\mathbf{R}}

\newcommand{\Div}{{\rm div\ }}

\newcommand{\Langle}{\left\langle}
\newcommand{\Rangle}{\right\rangle}

\newcommand{\vepsi}{\varepsilon}

\newcommand{\tr}{{\rm tr}}

\newtheorem{remark}{Remark}[section]

\begin{document}

\title[THE PHASE FIELD METHOD FOR MOVING INTERFACE PROBLEMS]{The Phase Field Method for Geometric Moving Interfaces 
	and Their Numerical Approximations}
%\title[]{The Phase Field Method for Geometric Moving Interfaces and Their Numerical Approximations}

\author{Qiang Du\dag}
%\address{Department of Applied Physics and Applied Mathematics
%		and the Data Science Institute, Columbia University, New York, NY 10027, U.S.A.}
%\email{qd2125@columbia.edu}
\thanks{\dag Department of Applied Physics and Applied Mathematics
	and the Data Science Institute, Columbia University, New York, NY 10027, U.S.A. (qd2125@columbia.edu). 
	The work of this author was partially supported  by the NSF grant DMS-1719699,
	and US Army Research Office MURI grant W911NF-15-1-0562.}
\author{Xiaobing Feng\ddag}
%\address{Department of Mathematics, The University of Tennessee, Knoxville, TN 37996, U.S.A. }
%\email{xfeng@math.utk.edu}
\thanks{\ddag Department of Mathematics, The University of Tennessee, Knoxville, TN 37996, U.S.A. (xfeng@math.utk.edu).
	The work of the first author was partially supported by the NSF grants: DMS-0410266 and DMS-1620168.}


%\date{}

\begin{abstract}
	This chapter surveys recent numerical advances in the phase field method  
	for geometric surface evolution and related geometric nonlinear partial differential 
	equations (PDEs). Instead of describing technical details of various numerical methods 
	and their analyses,  the chapter presents a holistic overview about the main ideas of phase field 
	modeling, its mathematical foundation, and relationships between the phase field formalism and other mathematical
	formalisms for geometric moving interface problems, as well as the current state-of-the-art of 
	numerical approximations of various phase field models with an emphasis on discussing 
	the main ideas of numerical analysis techniques. The chapter also reviews recent development
	on adaptive grid methods and various applications of the phase field modeling and their numerical 
	methods in materials science, fluid mechanics, biology and image science.
\end{abstract}

\maketitle

\tableofcontents

%%%%%%%%%%%%%%%%

%\input section1_v2r.tex

\section{Introduction} \label{sec-1}
%Topics to be covered include (not limited)
%
%\begin{enumerate}
%        \item Idea and definition of the phase field method for geometric moving interfaces
%        \item Main issues and challenges.
%        \item Historic timeline/milestones of the development of the method.
%        \item Outline of the chapter.
%\end{enumerate}

The idea of {\em the phase field} method could be traced back to Lord Rayleigh, Gibbs and Van der Waals and was 
used to describe material interfaces during phase transitions. It represents material interfaces as thin 
layers of finite thickness over which material properties vary smoothly. Such a thin layer 
of the width $O(\varepsilon)$ is 
often referred as a diffuse interface, and by the design the exact interface is guaranteed 
to be within the thin layer. In other words, this is amount to smear the exact sharp interface 
into a thin diffuse interface layer.  For this very reason, the phase field 
method is also known as {\em the diffuse interface method} in the literature. 

The phase field method was first introduced to model solid-liquid phase transition in which 
surface tension and non-equilibrium thermodynamics behavior become important at the interface. 
Computationally, the phase field method has features in common with the level set method 
\cite{osher_sethian88,evans_spruck91,chen_giga_goto91,sethian99,osher_fedkiw03} in that 
explicit tracking of the interface can be avoided in 
both mathematical formulation and in numerical computations, which has advantages 
in mesh generation and in capturing topological changes of the interface. The phase field method 
uses an auxiliary phase variable/function $u^\varepsilon$, also known as the order 
parameter, to indicate phases. The phase function assumes distinct values in the bulk 
phases away from the diffuse interface (or the interfacial region); the after-sought 
interface can be identified with an intermediate level set (e.g. the zero-level set) of the 
phase function $u^\varepsilon$.  It should be noted that although the level set method 
and the phase field method are intimately related, they differ fundamentally because the former 
tracks the exact sharp interface without introducing any diffuse layer. 

It is clear that the idea of the phase field method does not restrict to the material 
phase transition problems,  it is applicable to any moving interface (or free boundary) problems. 
Indeed, in the past forty years the phase field (or the diffuse interface) method has been developed 
into a major and general methodology for moving interface problems arising from  
astrophysics, biology, differential geometry, image processing, multiphase fluid mechanics, chemical and 
petroleum engineering,
materials phase transition and solidification. One common theme of these application problems is that 
interfacial energy %the surface tension 
plays an important role in each of these moving interface evolutions. Mathematically, 
 interfacial energies such as the  surface tension
 are characterized by the curvature(s) of the interface such as the mean and 
Gauss curvatures. It turns out that the mean and Gauss curvatures can be conveniently
expressed in term of the phase function, which makes the phase field method very effective for 
modeling the interfacial energetics, particularly the surface tension effect. 

Moving interface problems under the influence of the surface tension
 or other interfacial energy  belong to a larger class of so-called 
{\em geometric moving interface problems} in which the motion of the interfaces is driven by some  
curvature-dependent {\em geometric law}. The phase field formulations of such problems often give  
rise interesting and difficult geometric partial differential equations (PDEs), which is a main 
subject of this chapter. Among many geometric moving interface problems, the best known one perhaps 
is the {\em mean curvature flow} (MCF) whose governing geometric law is: $V_n(t)=H_{\Gamma(t)}$, 
where $V_n$ and $H_{\Gamma(t)}$ respectively stand for the normal velocity and the mean curvature of the moving 
interface $\Gamma(t)$ at time $t$. As the MCF is purely a geometric problem, it can be described by
various formulations including parametric, level set and phase field formulations. It is well-known that
the best known phase field formulation for the MCF is the Allen-Cahn equation: 
$u^\varepsilon_t -\Delta u^\varepsilon + \varepsilon^{-2} \bigl( (u^\varepsilon)^3-
u^\varepsilon \bigr)=0$, which is the simplest geometric PDE. Here the small positive 
parameter $\varepsilon$ measures the width of the diffuse interface layer. One of the most important  
theoretical results for the phase field method is to prove that the zero-level set 
$\Gamma^\varepsilon(t):=\{x\in \mathbf{R}^d;\, u^\varepsilon(x,t)=0\}$ evolves according to
the geometric law of the MCF and $\Gamma^\varepsilon(t)$ converges (as a set) to 
the exact MCF interface $\Gamma(t)$ as $\varepsilon\to 0^+$ \cite{evans_soner_souganidis92}. 

Looking at the PDE, it may not be clear why the Allen-Cahn equation is related to the MCF as its phase 
field model. This was indeed an open question for more than ten years. It took a considerable amount 
of effort by many researchers before the above convergence result was proved \cite{kohn_sternberg89,keller89,bronsard_kohn91,deMottoni_schatzman95,chen92,evans_soner_souganidis92}. 
Although rigorous proofs of the 
convergence of the phase field models/formulations of moving interface problems is generally difficult,  
there is a formal procedure which has been widely used for constructing/deriving phase field 
model/formulation for a given moving interface problem. This procedure is based on an {energetic 
approach} in which one first construct/postulate a Helmholtz free energy functional associated 
with the underlying moving interface problem, a desired phase field model/formulation then is obtained 
as the gradient flow for the energy functional in an appropriately chosen topology.  The Helmholtz free 
energy functional often can be written as the sum of two parts: the first part is called the bulk 
energy and the second part is the interfacial (or mixing) energy.  As its name indicates,
the interficial/mixing energy measures the energy stored in the interficial layer 
and can be defined and computed with the help of the phase function. As expected, it would depend 
on the width of the interfacial layer.  On the other hand,  the bulk energy is usually problem-dependent 
and may consist of several other types of energies such as kinetic and potential/gradient energy. 
It should be emphasized that this energetic approach makes the phase field modeling become 
systematic and methodical instead of ad hoc. It shifts the main task of modeling to the construction of the Helmholtz free energy functional and  the driving forces that provide the mechanism of energy dissipation.

After the phase field model is obtained, several important questions naturally arise and must be addressed. Among
them we mention the following three:  (1) Is the model mathematically well-posed as a geometric PDE problem? 
(2) Does the phase field model converge to the original sharp interface model as the width of the interfacial layer
goes to zero? If it does, how fast or what is the convergence rate? (3) How to efficiently solve the phase field 
model (i.e., the geometric PDE problem)  numerically?  It turns out that these questions are not all
easy to answer. For  question (2) the main difficulties are to derive required uniform in $\varepsilon$ 
estimates for the phase function and to characterize the solution of the limiting 
model. Such a rigorous proof of convergence is still missing for many moving interface 
problems from materials science and multiphase fluid mechanics. 
For question (3) several  important issues must be considered. First, in order to preserve the 
gradient flow structure at the discrete level, it is desirable to ensure the numerical method to be energy stable 
besides to be accurate, which is a non-trivial task, especially if one aims to  have high order
(time) schemes (see \cite{yang16,yang2017numerical,shen2018} for some recent advances  
in this direction). Second, the phase field method is used to formulate/approximate a moving 
interface problem (i.e. the underlying sharp interface problem) and numerical methods are  
constructed to solve the phase field model. If the goal is to solve the sharp interface problem,
a critical question is that whether the numerical solution of the phase field model
converges to the solution of the underlying sharp interface model as the width 
of the diffuse interface tends to zero (i.e., $\varepsilon\to 0$). One may further ask in what
sense the above convergence should be interpreted. It turns out that these
issues are quite difficult to address and remain open for many phase field models and their numerical
approximations. The successful proofs all require to carry out some delicate and nonstandard 
error analyses \cite{Feng_Prohl03b,Feng_Prohl04a,Feng_Prohl04c,Feng_Prohl05a,Feng_Li15,FLX16a}. 
Finally, in order to resolve the thin interfacial layer, one must use spatial (and temporal) 
mesh sizes that are much smaller than the layer width $\varepsilon$ for computer simulations 
(note that $\varepsilon\to 0$). For such extremely small mesh sizes, the resulting large linear and nonlinear 
algebraic problems require huge computational effort to solve. This is especially 
true in the three-dimensional case with uniform spatial meshes. To overcome the difficulty, one may
rely on highly scalable algorithms  \cite{GB2011,zhang16}  or use adaptive 
 methods that only use fine meshes in the interfacial layer and much coarser meshes 
away from the layer. The use of adaptive mesh is a natural choice considering the fact that 
the phase function has a nice structured profile that assumes distinct (close-to-constant) values in bulk phases 
away from the diffuse interface and has a large gradient in the diffuse interface. 

Although the phase field method was introduced much earlier as a modeling tool for material phase
transition (in particular, solidification),  it was proposed as a computational 
technique in the early 1980s by J. S. Langer, G. Fix, G. Caginalp and others. 
G. Fix \cite{fix83} was perhaps the first to use the phase field method to 
numerically solve moving interface (or free boundary) problems. 
Since then the phase field method has garnered a lot of interest and become more and more 
popular. It has been developed into a general numerical methodology for solving various moving 
interface (i.e., free boundary) problems from differential geometry, materials science, fluid mechanics, 
and biology, just to name a few. Practically, the methodology is well understood,  however,  
theoretically,  many difficult questions remain to be addressed both at the PDE level and
at the discrete level. Moreover, computer simulations of the phase field method, 
especially  in high dimensions, remains to be challenging because of shear amount of computations 
involved. Indeed, phase field simulations on extreme scales have been used
as a test of computing power on world's leading high performance computers as exemplified by the 
appearances on the Gordon
Bell prize competition \cite{GB2011} (2011 Gordon Bell prize winner) 
and \cite{zhang16} (2016 Gordon Bell prize finalist, which reported
the largest 3D phase field simulations to date that 
reached 40 Petaflops on the world's fastest supercomputer at the time, employing
over 300 billion spatial grid points ($6700^3$) and utilizing over 10 million cores). With the 
increasing popularity of phase field modeling
in more and more application domains, it is imperative that
efficient numerical algorithms and solvers such as adaptive grid methods and 
preconditioning techniques must be developed and utilized.

The phase field method provides a diffusive interface description to many 
physical processes and geometric models. Its capability goes beyond co-dimensional one 
surfaces. For example, using vector fields or tensor fields, they can be used to describe the
characteristics of geometric information with higher co-dimension such as point defects in
two dimensions and curves in three dimensions. Mathematical models of this kind include the
celebrated Ginzburg-Landau models of superconductivity, Landau-deGennes model for liquid crystal
flow, generalized Ohta-Kawasaki model of diblock and triblock copolymers,  among others
\cite{davis98,du92,kim05,ohta86,priestly12,tinkham04}.  While many numerical methods designed for scalar
equations can be extended to systems, there are new challenges associated with 
invariance and symmetries and choices of gauge \cite{du05} in addition to added computational
complexity due to additional quantities of interests.

This chapter consists of additional seven sections. Section \ref{sec-2} devotes to the 
mathematical foundation of the phase field method. This includes its formulations for specific 
geometric moving interface problems, its intimate connection to the level set method, convergence 
of the phase field method as the defuse interface width tends to zero, fine properties 
of solutions of phase field models, and phase function representations of geometric quantities. 
Section \ref{sec-3} focuses on time-stepping schemes for the phase field models, in addition
to ensure their accuracy, the emphasis is on presenting various energy stable schemes 
of different orders which include the classical energy splitting schemes and stabilized schemes
as well as the newly developed IEQ (Invariant Energy Quadratization) and SAV (Scalar Auxiliary Variable)
schemes. 
Section \ref{sec-4} devotes to spatial discretization methods, all the classical methodologies are 
covered.  Among them are finite difference methods, finite element methods, spectral methods,
discontinuous Galerkin methods, and isogeometric methods.  
Section \ref{sec-5} focuses on two convergence theories of fully discrete numerical methods for 
phase field methods.  The first theory analyzes rates of convergence (or a priori error estimates)
for a fixed diffuse interface width $\varepsilon$, while the second theory considers the convergence 
to the solution of the underlying sharp interface problem when both the mesh parameters and the 
physical parameter $\varepsilon$ all tends zero. 
Section \ref{sec-6} addresses a posteriori error estimates which are used to design adaptive methods.  
The emphasis will be on the residual-based a posteriori estimators for finite element discretizations. 
All three types of popular adaptive methodologies will be discussed, they include $h$-, $hp$- and $r$-adaptive 
methods. 
Section \ref{sec-7} lists a host of applications of the phase field method ranging  from materials science, 
fluid mechanics to biology. 
It also discusses some extensions of the classical phase field method such as nonlocal 
and stochastic phase field models. 
The chapter is ended with a brief summary and a few concluding remarks given in Section \ref{sec-9}.

\section{Mathematical Foundation of the Phase Field Method} \label{sec-2}

%Topics to be covered include (not limited)
%\begin{enumerate}
%        \item Geometric laws and geometric moving interfaces.
%        \item Phase field formulations/models of some well-known geometric laws.
%        \item Convergence of the phase field method (for the mean curvature flow).
%        \item Fine properties of solutions of phase field models.
%        \item Phase function representations of geometric quantities.
%\end{enumerate}

%\bigskip
%\hrulefill

\subsection{Geometric Surface Evolution}
Mathematically, a hypersurface or an interface in the Euclidean space $\bR^{d+1}$ refers to 
co-dimension one subset $\Gamma$ of $\bR^{d+1}$.  A moving interface (or a free boundary)
refers to one-parameter family of interfaces $\{\Gamma_t\}_{t\geq 0}$ which evolves in the 
space $\bR^{d+1}$ and are governed by some explicit or implicit mechanism/law, often called a 
{\em geometric law}. The parameter $t$ represents the time. Such a geometric law determines 
how the interface moves. Most commonly seen geometric laws specify the velocity $V$ or 
normal velocity $V_n$  of the points (thinking them as particles) on the interface at each given 
time $t$.  It may depend
on intrinsic features (such as curvatures) of the  interface and/or on external factors 
(such as flow velocity) of the environment where the interface exists.  

Specifically, geometric surface evolution concerns with the following question: Given an
initial hypersurface $\Gamma_0\in \bR^{d+1}$, find $\Gamma_t\in \bR^{d+1}$	such that
\begin{equation}\label{eq2.1}
V_n(t)= F_{\mbox{\tiny{int}}} (\lambda_1,\cdots,
\lambda_d) + F_{\mbox{\tiny{ext}}}\qquad\mbox{on }\G_t ,
\end{equation} 
where $V_n=V\cdot n$ denotes the (outward) normal velocity of $\Gamma_t$ and $n=n_{\G_t}$ stands for the 
outward normal to $\Gamma_t$. $\{\lambda_j(t)\}_{j=1}^d$ are the principle curvatures of $\Gamma_t$ and 
$F_{\mbox{\tiny{int}}}$ is a given function of $\{\lambda_j\}_{j\geq 1}$.  $F{\mbox{\tiny{ext}}}$
denotes an external (source) function. The geometric law \eqref{eq2.1} says that the normal velocity 
$V_n$ of the interface is determined/driven by the sum 
of an internal ($F_{\mbox{\tiny{int}}}$) and an external ($F{\mbox{\tiny{ext}}}$) ``forces". 
We note that problem \eqref{eq2.1} is stated purely as a geometric problem, it may be  
embedded in a more complicated moving interface problem arisen from a particular application.

\subsection{Examples of Geometric Surface Evolution}\label{sec-1.2.2}
The following is a list of some best known geometric moving interface problems in the literature. 
We rephrase them in the framework of problem \eqref{eq2.1} by explicitly describing the functions 
$F_{\mbox{\tiny{int}}}$ and $F{\mbox{\tiny{ext}}}$ for each example.

{\bf Example 1:} ({\em Mean curvature flow})  The mean curvature flow (MCF) is defined 
	by setting 
	\[
	F_{\mbox{\tiny{int}}} = H:=\sum_{j=1}^d \lambda_j ,
	\qquad F_{\mbox{\tiny{ext}}}\equiv 0    
	\]
	in \eqref{eq2.1}.  Where $H=H(t)$ denotes the mean curvature of $\Gamma_t$. So the MCF 
	seeks a family of hypersurfaces whose normal velocity at every point on the hypersurface  
	is equal to the mean curvature of the hypersurface at that point for all time $t>0$. 
	The MCF is a curve-shortening (in 2-D) and an area minimization flow (in higher dimensions)
	in the sense that it can be interpreted as a gradient flow for the curve length functional (in 2-D) 
	and the surface area functional (in higher dimensions) \cite{zhu02,ecker04,giga06}. 
	The MCF finds applications 
	in many fields such as materials science, image processing, and multiphase fluids. 
	
\bigskip
{\bf Example 2:} ({\em A generalized mean curvature flow})  The generalized mean curvature 
	flow (gMCF) to be considered in this chapter refers to the case with 
	\[
	F_{\mbox{\tiny{int}}} = H, %:=\sum_{j=1}^d \lambda_j ,
	\qquad F_{\mbox{\tiny{ext}}}\equiv \bv\cdot n_{\G_t} +g 
	\]
	in \eqref{eq2.1}. Clearly, the difference between and the MCF and the gMCF is that 
	the latter has a non-zero external force $F_{\mbox{\tiny{ext}}}$. Here $\bv$ 
	represents the background fluid velocity, and $g$ denotes the combined other external forces. 
	
\bigskip 
{\bf Example 3:}  ({\em Inverse mean curvature flow}) The inverse mean curvature flow (IMCF) 
     is defined by setting  
	\[
	F_{\mbox{\tiny{int}}}= \frac{1}{H},\qquad F_{\mbox{\tiny{ext}}}=0
	\]
	in \eqref{eq2.1}, it describes the evolution of a family of hypersurfaces whose normal velocity at 
	every point on the hypersurface is equal to the reciprocal of its mean curvature  
	at that point for all time $t>0$. The IMCF finds application in general relativity 
	and was used as a main tool to prove the {\em Penrose Inequality} by Huisken and Ilmanen 
	\cite{huisken_ilmanen01}.
	
\bigskip
{\bf Example 4:} ({\em Surface diffusion flow}) The surface diffusion flow refers to the
     case with 
	\[
	F_{\mbox{\tiny{int}}} =\Delta_{\G_t} H,\qquad F_{\mbox{\tiny{ext}}}=0
	\]
	in \eqref{eq2.1}. Where $\Delta_{\G_t}$ denotes the surface Laplace operator on $\G_t$, 
	so the flow requires that the normal velocity of the hypersurface equals the surface Laplace  
	of its mean curvature at every point on the surface for all time $t>0$. The surface diffusion 
	arises from applications in materials science, image processing, and cell biology \cite{chopp_sethian99,
	escher_mayer_simonett98,elliott2010surface}.   
	
\bigskip
{\bf Example 5:} ({\em Willmore flow}) The Willmore flow is defined by setting 
	\[
	F_{\mbox{\tiny{int}}}=-\Delta_{\G_t} H-2H(H^2-K),
	\quad K:=\prod_{j=1}^d \lambda_j,
	\quad F_{\mbox{\tiny{ext}}}=0
	\]
	in \eqref{eq2.1}. Note that $K$ denotes the Gauss curvature of the surface.  The Willmore flow 
	can be interpreted as a gradient flow for the total quadratic mean curvature functional $\int_{\G_t} H^2$ 
	and find applications in biology and materials science \cite{simonnett01,du2004phase,wang2008modelling,zhang2009phase, elliott2010surface, bretin2015phase,marques_neves14}.
	 
\bigskip
{\bf Example 6:} ({\em Hele-Shaw flow}) To define the Hele-Shaw flow, we set 
	\[
	F_{\mbox{\tiny{int}}}=0, \qquad F{\mbox{\tiny{ext}}}
	= \frac12 \left[\frac{\p w}{\p n}\right]_{\G_t}, 
	\]
	in \eqref{eq2.1}, where $w$ is defined by 
	\begin{alignat*}{2}
	\Del w&=0 &&\quad\mbox{in } \Ome\setminus {\G_t}, \\
	w&=\sigma H &&\quad\mbox{on } {\G_t},
	\end{alignat*}
	here $\Omega\in \mathbf{R}^d$ is a bounded domain, $\frac{\p w}{\p n}$ stands for the normal derivative of $w$ and $[w]_{\G_t}$ 
	denotes the jump of a function $w$ across the interface $\G_t$.  $\sigma$ is a positive constant and is 
	called the surface tension in the literature. Note that we omit the outer boundary condition for $w$. 
	The Hele-Shaw flow says that the normal velocity of the 
	interface equals the jump of the normal derivative of the pressure field $w$ across the 
	interface.  Since the pressure depends on the mean curvature $H$, then the normal velocity implicitly 
	depends on the mean curvature $H$.  The Hele-Shaw flow arises in the study of two-phase fluids. 
	It should be noted that the identical model also arises in materials sciences and is known as the 
	Mullins-Sekerka model, though function $w$ stands for the temperature, not the pressure, in the 
	Mullins-Sekerka model \cite{hele-shaw1898b,mullins_sekerka63}.
 
\bigskip
{\bf Example 7:} ({\em Generalized Stefan problem}) The generalized Stefan problem can be defined by setting 
	\[
	F_{\mbox{\tiny{int}}}=0, \qquad 
	F_{\mbox{\tiny{ext}}} =\frac12 \left[\frac{\p \phi}{\p n}\right]_{\G_t}, 
	\]
	in \eqref{eq2.1}, where $\phi$ is defined by 
	\begin{alignat*}{2}
	\frac{\p \phi}{\p t} &-\Delta \phi =0 &&\quad\mbox{in } \Ome\setminus \G_t, \\
	\qquad \phi &=\sigma (H - \alpha V_n) &&\quad\mbox{on } {\G_t},
	\end{alignat*}
	here $\alpha$ is another positive constant. Compared with the Hele-Shaw flow, the temperature $\phi$ now 
	satisfies the heat equation in $\Omega\setminus \G_t$ and the normal velocity $V_n$ appears 
	in the boundary condition for $\phi$ (thus, it appears in both sides of \eqref{eq2.1}). Note that 
	we omit the outer boundary condition and the initial condition for $\phi$.  The generalized 
	Stefan problem describes the solidification process of fluids with consideration of the surface diffusion 
	(characterized by $\sigma$) and the super/under cooling effect (characterized by $\alpha$) 
	\cite{caginalp_chen98,langer86,fix83}. 
	
\bigskip
{\bf Example 8:} ({\em Two phase immiscible fluids}) The sharp interface two phase immiscible 
flow model is often written as 
	%\begin{align*}
	%&\mbox{momentum equations in bulk of each fluid phase, NS-equations or non-Newtonian 
	%fluid equations}, \\
	%&\mbox{interface conditions: (dis)continuity of normal stress and velocity}.
	%\end{align*}
	\begin{alignat*}{2}
	\rho \Bigl( \frac{\p \bu}{\p t}+ (\bu\cdot\nab)\bu \Bigr)
	-\Div T(\bu,p)&=\mathbf{f} &&\qquad\mbox{in }\Ome_T\setminus \G_t,\\
	\Div \bu &=0  &&\qquad\mbox{in }\Ome_T\setminus \G_t,\\
	[T(\bu,p) n]_{\G_t} &= \sigma H &&\qquad\mbox{on } \G_t, \\
	[\bu]_{\G_t} & = 0  &&\qquad\mbox{on } \G_t,\\
	%+\mbox{BC and IC} & &&
	T(\bu,p):= \nu \bigl( \nab\bu+\nab\bu^T \bigr) &-  p I,  &&
	\end{alignat*}
where $\Omega_T:=\Omega\times (0,T)$, $\bu$ and $p$ denote the fluid velocity and pressure. $\nu (>0)$ 
is the viscosity coefficient. Again, $[\cdot]_{\G_t}$ denotes the jump function. 
The first two equations in the PDE system are called the Navier-Stokes equations. $T(\bu,p)$ is called the 
deformation tensor which can be replaced by more complicated nonlinear version in the
case when complex fluids are involved \cite{jacqmin99,anderson_mcfadden_wheeler98,lowengrub1998quasi,liu_shen03,yue2004diffuse,feng2005energetic}. 
To fit the above model into the form \eqref{eq2.1}, we simply set 
\[
F_{\mbox{\tiny{int}}}=0, \qquad 
F_{\mbox{\tiny{ext}}} =\bigl[ T(\bu, p)  n\bigr]_{\G_t}.
\]
It follows from \eqref{eq2.1} and the first interface  condition in the PDE system is 
$V_n=\sigma H$,  hence the fluid interface $\G_t$ evolves as a scaled mean curvature flow.  
 
\bigskip
{\bf Example 9:} ({\em Two phase Hele-Shaw flows}) An interesting special case of two phase flows is 
so-called the two phase Hele-Shaw flow which refers to the motion of (one or more) viscous fluids between two flat 
parallel plates separated by an infinitesimally small gap. Such a physical setup is often called a Hele-Shaw cell
and was originally designed by Hele-Shaw to study two dimensional potential flows\cite{hele-shaw1898b,kim05, Feng_Wise12,wang_zhang13,han_wang18}. The sharp interface model of two phase Hele-Shaw flows is given by
\begin{alignat*}{2}
\bu &=-\frac{1}{12 \eta}\bigl(  \nabla p -\rho \mathbf{g} \bigr)  &&\qquad\mbox{in }\Ome_T\setminus \G_t,\\
\Div \bu &= 0 &&\qquad\mbox{in }\Ome_T\setminus \G_t, \\
[p] &= \gamma H &&\qquad\mbox{on } \G_t, \\
[\bu\cdot n ] & =0 &&\qquad\mbox{on }  \G_t,
\end{alignat*}
where the first interface condition is called the Laplace-Young condition in which $\gamma$ is the dimensionless 
surface tension coefficient and $H$ stands for the mean curvature of $\G_t$. To fit the above model into the 
form \eqref{eq2.1}, we set 
\[
F_{\mbox{\tiny{int}}}=0, \qquad 
F_{\mbox{\tiny{ext}}} =[ p]_{\G_t},
\]
which says that the fluid interface $\G_t$ evolves as a scaled mean curvature flow.

%%%%%%%
\subsection{Mathematical Formulations and Methodologies}
As equation \eqref{eq2.1} is a geometric problem, in order to analyze and approximate its solutions, 
one needs to choose an appropriate formulation and to define the meaning of solutions, which in turn 
then lead to different mathematical and numerical methods for the problem. 
	 
Several  mathematical formulations of \eqref{eq2.1} have been proposed and developed 
in the past forty years.  Among them the best known one  
is the parametric formulation (cf. \cite{ecker04, zhu02} and the references therein), 
in which the coordinates of all points on the hypersurface are written 
as functions of chosen parameters. All best know geometric quantities can be computed in terms 
of the coordinate functions. The parametric formulation has been widely used 
to study smooth geometric surface evolutions. It has also been used to develop numerical 
methods for computing the solutions (i.e., the coordinate functions).  
In order to deal with non-smooth geometric surface evolutions,  other formulations and weak solution 
concepts have been introduced and investigated. One of such formulations is the varifold formulation 
developed by Brakke \cite{Brakke78}, especially for the mean curvature flow.  Its measure-valued 
weak solutions are called varifold solutions. Another weak formulation and weak solution 
notion, called the theory of minimal barriers, were introduced by De Giorgi \cite{degiorgi94} and were 
further developed by Bellettini and Novaga \cite{bellettini_novaga97, bellettini13} including establishing the
connection between barrier solutions and other type weak solutions.  
The fourth formulation, which has been very popular for both PDE analysis and numerical approximation, 
is the level set formulation \cite{osher_sethian88,osher_fedkiw03,sethian99}. 
The corresponding weak solutions are called level 
set solutions \cite{evans_spruck91,chen_giga_goto91}.
The level set formulation provides a convenient and effective formalism not only for analyzing 
curvature-driven flows such as the mean curvature flow but also for approximating their (level set) 
solutions numerically (see more details below).
Finally, the fifth formulation of \eqref{eq2.1} is the phase field formulation \cite{fix83,langer86, fife88, mcfadden02,
anderson_mcfadden_wheeler98}, which is the 
main subject of this chapter. Like the level set formulation, the phase field formulation
also provides a convenient and effective formalism for both mathematical analysis and numerical
approximation of problem \eqref{eq2.1}, especially, in the case when $F_{\mbox{\tiny{ext}}}\not\equiv 0$.  
 
%%%%%%%
\subsection{Level Set and Phase Field Formulations of the MCF}
We now use the mean curvature flow (MCF) as an example to demonstrate the derivations of both 
the level set and phase field formulations of the MCF. 
 
The level set method was introduced by Osher and Sethian in \cite{osher_sethian88} for 
general moving interface problems. 
As its name indicates, the main idea of the level set formulation/method is to represent the 
hypersurface $\G_t$ as a zero-level set of a function $u$ in $\mathbf{R}^{d+1}$, that is 
\begin{equation}\label{level_set}
	\G_t:=\{x\in\Ome;\, u(x,t)=0\}, 
\end{equation}
and then to evolve the level-set function $u$, instead of the interface $\G_t$.  To transfer 
equation \eqref{eq2.1} for $\G_t$ to an equation for $u$, we formally differentiate 
equation $u(x,t)=0$ with respect to $t$, while treating $x=x(t)$ as an implicit function, 
and using the chain rule to get 
\begin{align*}
\frac{\p u}{\p t} +\nabla u\cdot\frac{dx}{dt}=0.
\end{align*}
Since $V = \frac{dx}{dt}$ is the velocity of the surface, then
\begin{align}\label{level_set_1}
\frac{\p u}{\p t} +\nabla u\cdot V=0.
\end{align}
Equation \eqref{level_set_1} is often called the level set equation, and it is determined by the 
velocity field $V$ and the initial condition $u_0$ such that $\Gamma_0 = \{x\in \mathbb{R}^{d+1}; u_0(x)=0\}$.
To illustrate, we consider the mean curvature flow described in {\bf Example 1} whose geometric law is
\begin{align*}
V_n(t, \cdot) = H(t, \cdot). 
\end{align*}
By the differential geometry facts we have
\begin{align*}
n = -\frac{\nabla u}{|\nabla u|} \qquad\text{and}\quad H = -\mbox{div}(n),
\end{align*}
where $n$ stands for the inward normal to $\G_t$. Then the level set equation \eqref{level_set_1} becomes
\begin{align*}
0=\frac{\p u}{\p t}+\nabla u\cdot V=\frac{\p u}{\p t}-|\nabla u|V_n
=\frac{\p u}{\p t}-|\nabla u|H
=\frac{\p u}{\p t} - |\nabla u|\,\mbox{div}\left(\frac{\nabla u}{|\nabla u|}\right),
\end{align*}
or
\begin{align}\label{level_set_2}
\frac{\p u}{\p t} - |\nabla u|\, \mbox{div} \left(\frac{\nabla u}{|\nabla u|}\right)=0.
\end{align}
Equation \eqref{level_set_2} is the well-known level set formulation of the mean curvature 
flow \cite{evans_spruck91,chen_giga_goto91,giga06}. 

It should be pointed out that here we only showed that equation \eqref{level_set_2} holds on $\G_t$, but it easy 
to verify that it actually holds on all level sets, hence, holds in $\bR^{d+1}$. It is also 
important to note that the level-set equation \eqref{level_set_1} remains the same for {\em all} moving 
interface problems although the velocity $V$ may be different for different problems.  

The phase field method for the moving interface problems can be traced back to 
Lord Rayleigh \cite{Rayleigh1892}, Van der Waals \cite{waals1892} and Gibbs \cite{EBW91}. It was originally 
developed as a model for solidification, but has been used for many other applications, 
such as crack propagation, electromigration, crystal and tumor growth. The main idea of phase 
field method is to introduce thickness to the interface, more precisely, it seeks a phase field 
function $u^{\eps}$ such that the interface lies in the narrow region (called the diffuse interface)
\begin{align}\label{phase_field}
\Gamma_t\subset Q_t^{\eps} := \bigl\{x(t)\in\mathbf{R}^{d+1}: |u^{\eps}(x(t),t)|\leq1-\mathcal{O}(\varepsilon)\bigr\}.
\end{align}
Here $\varepsilon$ is a small positive constant, which controls the width of $Q_t^{\varepsilon}$. 
The phase field function takes two distinct value $+1$ and $-1$, which represent two distinct phases, 
with a smooth change between $-1$ and $+1$ in $Q_t^{\eps}$. The zero level set 
$\Gamma_t^{\varepsilon} := \{x(t)\in\mathbf{R}^{d+1}; u^{\varepsilon}(x(t),t)=0\}$ of $u^{\varepsilon}$, 
which is contained in the diffuse interface $Q_t^{\varepsilon}$, is often chosen to represent 
$\Gamma_t$ approximately. Like the level set method, this diffuse interface approach provides 
a convenient mathematical formalism for numerically approximating moving interface problems 
because explicitly tracking the interface is not needed in the formulation.
Unlike the level set method, there is no master phase field equation which is valid for all 
moving interface problems, instead, the phase field formulation is problem-dependent and are often  
difficult to derive. The difficulty is due to the fact that the interface lies inside 
$Q_t^{\varepsilon}$, but the specified location of the interface is unknown, so the curvature at 
the interface can not be calculated exactly as in the level set method. Below we again use the mean 
curvature flow (MCF) as an example to show a formal derivation of its phase field formulation when
the profile of the interface is postulated.   
To this end, let $d(x)$ denote the signed distance function between point $x$ and the interface 
$\Gamma_t$, and consider the fact that the solution approximates the $\tanh(\cdot)$ function, 
we heuristically postulate that  
\begin{align}\label{phase_field_1}
u^{\eps}(x,t) := \tanh\bigg(\frac{d(x)}{\sqrt{2}\eps}\bigg),
\end{align}
because it matches with the desired profile for the phase function $u^\eps$. Then we have
\[
\tanh^{\prime}(s)=1-\tanh^2(s),\quad
\tanh^{\prime\prime}(s)=-2\tanh(s)(1-\tanh^2(s)),
\]
and
\begin{align*}
\nabla u^{\eps}(x) &= \frac{\tanh\bigl(\frac{d(x)}{\sqrt{2}\eps}\bigr)}{\sqrt{2}\eps}\nabla d(x), \\
D^2 d(x) &= \frac{\sqrt{2}\eps}{1-(u^{\eps}(x))^2}\bigg( D^2 u^{\eps}(x)+\frac{2u^{\eps}(x)}{1-(u^{\eps}(x))^2}\nabla u^{\eps}(x)\otimes\nabla u^{\eps}(x)\bigg). 
\end{align*}
By differential geometry facts we have
\[
|\nabla d(x)| = 1, \quad 
|\nabla u^{\eps}(x)|^2 = \frac{1}{2\eps^2}\Bigl(1-(u^{\eps}(x))^2 \Bigr)^2, 
\]
hence,
\begin{align}\label{phase_field_2}
H = \tr(D^2 d(x)) = \frac{\sqrt{2}\eps}{1-(u^{\eps}(x))^2}\bigg(\Delta u^{\eps}(x)+\frac{1}{\eps^2}(u^{\eps}(x)-(u^{\eps}(x))^3)\bigg).
\end{align}

Recall that the approximate moving interface $\G_t^\eps$ is represented by the zero-level set of
$u^\eps_t$, that is, $u^\eps_t(x,t)=0$ on $\G_t^\eps$. As in the derivation of the level set 
equation, we formally differentiate equation $u^\eps(x,t)=0$ with respect to $t$, while treating 
$x=x(t)$ as an implicit function, and using the chain rule to get 
\begin{align}\label{phase_field_3}
0=\frac{\p u^\eps}{\p t}+\nabla u^\eps\cdot V=\frac{\p u^\eps}{\p t}-|\nabla u^\eps|V_n
=\frac{\p u^\eps}{\p t}-|\nabla u^\eps|H,
\end{align}
here we have used the facts that $V=\frac{d x(t)}{d t}$ and $n=-\frac{\nabla u^\eps}{|\nabla u^\eps|}$. 
Combining \eqref{phase_field_2} and \eqref{phase_field_3}, we obtain the follow phase field equation
for the MCF: 
\begin{align}\label{phase_field_4}
\frac{\p u^{\eps}}{\p t} - \Delta u^{\eps} + \frac{1}{\eps^2}\bigl((u^{\eps})^3-u^{\eps}\bigr)=0.
\end{align}
%It was proved in \cite{Evans_Soner_Souganidis92} that $\Gamma_t^{\epsilon}$ converges 
%to $\Gamma_t$ defined in \eqref{level_set_2} as $\epsilon\rightarrow 0$.

Equation \eqref{phase_field_4} is called the Allen-Cahn equation in the literature, it was 
introduced by Samuel M. Allen and John W. Cahn in \cite{allen79}  
as a model to describe the phase separation process of a binary alloy at a fixed temperature. 
In the original Allen-Cahn equation $u^\eps$ denotes the concentration of one of 
the two species of the alloy. We also remark that equation \eqref{phase_field_4} differs from the
original Allen-Cahn equation in the time scale, $t$ here represents $\frac{t}{\eps^2}$ 
in the original Allen-Cahn equation, hence, it is a fast time. 
%The Allen-Cahn equation is not mass-conserved because $\int_\Ome u\, dx$ is not a constant in $t$.

The Allen-Cahn equation can also be derived using the modern {\em energetic approach}.  To
this end, we introduce the following Cahn-Hilliard free-energy functional
\begin{equation}\label{CH_functional} 
\mathcal{J}(u):= \int_\Omega \Bigl( \frac12 |\nabla u|^2 + \frac{1}{\eps^2} F(u) \Bigr) \, dx, \qquad
F(u):=\frac14 \bigl( u^2-1 \bigr)^2. 
\end{equation}
where the first term in $\mathcal{J}$ is called the bulk energy and the second term is called 
the interficial (or potential) energy. Then \eqref{phase_field_4} can be interpreted as 
the $L^2$-gradient flow for $\mathcal{J}$, that is,
\begin{align}\label{phase_field_5}
\frac{\p u}{\p t} = -\mathcal{J}^\prime (u) \qquad\mbox{in } L^2(\Omega), 
\end{align}
where $\mathcal{J}^\prime (u)$ denotes the G\^ateaux derivative of $\mathcal{J}$ at $u$ 
in the specified topological space.  For 
the $L^2(\Omega)$ space, it is easy to verify that 
\begin{align}\label{phase_field_6}
\mathcal{J}^\prime (u) = -\Delta u + \frac{1}{\eps^2} f(u), \qquad 
\end{align}
where $f(u):=F'(u)=u^3-u$ for the double well potential given in \eqref{CH_functional}.
 We thus get a concrete form of \eqref{phase_field_5} given by
\begin{align}\label{phase_field_5a}
\frac{\p u}{\p t} = \Delta u - \frac{1}{\eps^2} f(u). \qquad 
\end{align}

Clearly,  the above energetic approach is quite simple and systematic compared to the traditional 
derivation which could be involved for complex moving interface problems. The key ingredient/step 
is to design a ``correct" energy functional $\mathcal{J}$  and to choose the ``right" topology 
(or inner product) for a specific problem,  the resulting phase field equation is then given by the 
general gradient flow equation \eqref{phase_field_5}.  For this reason we may regard equation  
\eqref{phase_field_5} as the (master) phase field equation. 

%%%%%%%
\subsection{Phase Field Formulations of Other Moving Interface Problems}
In this subsection we list the known phase field formulations for all moving interface 
problems introduced in Section 1.2.2.  

For the generalized mean curvature flow defined in {\bf Example 2}, its phase field formulation 
is given by the following {\em convective} Allen-Cahn equation \cite{caginalp_chen98,chen1998convergence}:
\begin{equation}\label{convective_AC}
	\frac{\p u^\eps}{\p t} - \Del u^\eps 
	+ \mathbf{v}\cdot \nab u^\eps 	+ \frac{1}{\eps^2} f(u^\eps)=g.
\end{equation}
 
For the surface diffusion flow defined in {\bf Example 4}, its phase field 
formulation is the following {\em degenerate} Cahn-Hilliard equation \cite{cahn_elliott_cohn96}:
\begin{align}\label{degenerate_CH}
\frac{\p u^\vepsi}{\p t}+\Div\Bigl(b(u^\vepsi)\nab \bigl(\Del u^\vepsi-\frac{1}{\vepsi^2}f(u^\vepsi) \bigr) \Bigr)=0,
%\mu^\vepsi &= -\Del u^\vepsi +\frac{1}{\vepsi^2} f(u^\vepsi),
\quad b(z):=1-z^2
\end{align}
where $b(z):=1-z^2$. 

For the Willmore flow introduced in {\bf Example 5},  its phase field formulation 
is given by the following fourth order PDE \cite{du2004phase, bellettini13, bellettini_novaga97, BMO15}: 
\begin{align}\label{Willmore_PF}
\frac{\p u^\vepsi}{\p t} +\Del \Bigl( \Del u^\vepsi -\frac{1}{\vepsi^2} f(u^\vepsi) \Bigr) 
-\frac{1}{\eps^2} f'(u)  \Bigl( \Del u^\vepsi -\frac{1}{\vepsi^2} f(u^\vepsi) \Bigr)  =0,
\end{align}
corresponding to the following phase field relaxation of the Willmore energy functional:
\begin{equation}\label{Willmore_energy}
\mathcal{W}(u):= \frac{1}{2\varepsilon} \int_{\Omega}  \Bigl(\varepsilon\Del u^\vepsi -\frac{1}{\vepsi} f(u^\vepsi)  \Bigr)^2\, dx. 
\end{equation}

For the Hele-Shaw flow defined in {\bf Example 6}, its phase field formulation is
given by the following well-known Cahn-Hilliard equation \cite{cahn58, Pego89, alikakos94}: 
\begin{equation}\label{cahn_hilliard}
\frac{\p u^\vepsi}{\p t} +\Del \Bigl( \vepsi \Del u^\vepsi -\frac{1}{\vepsi} f(u^\vepsi) \Bigr) =0,
\end{equation}
which can also be obtained from \eqref{degenerate_CH} after setting $b(z)\equiv 1$. 
 
For the generalized Stefan problem introduced in {\bf Example 7} of Section \ref{sec-1.2.2},  
its phase field formulation is given by the following so-called classical phase field model
\cite{caginalp_chen98}:
\begin{align}\label{classical_PF} 
\alpha(\vepsi)\frac{\p u^\vepsi}{\p t}&= \vepsi\Del u^\vepsi
-\frac{1}{\vepsi} f(u^\vepsi) + s(\vepsi) \phi^\vepsi \\
c(\vepsi)\frac{\p \phi^\vepsi}{\p t}&= \Del \phi^\vepsi-\frac{\p u^\vepsi}{\p t}.
\end{align}
Note that $u^\varepsilon$ denotes the phase function in the model. 
 
	\begin{remark}
		There is also a phase field model for generalized anisotropic Stefan problem which is widely used in the materials 
		science community \cite{provatas2011phase, mcfadden02, boettinger2002phase, plapp_karma00}. 
	\end{remark}
 
For the two-phase-fluid diffuse interface model introduced in {\bf Example 8} 
of Section \ref{sec-1.2.2}, when $\rho=1$ and $\nu=1$, its phase field formulation is given by 
the following coupled Navier-Stokes and Cahn-Hilliard system \cite{anderson_mcfadden_wheeler98, jacqmin99,liu_shen03}:
\begin{alignat}{2} \label{NS_CH} 
\frac{\p \bu^\varepsilon}{\p t} + (\bu^\varepsilon\cdot\nab)\bu^\varepsilon
-\Delta \bu^\varepsilon +\nabla p^\varepsilon &=\mu^\varepsilon\nabla \varphi^\varepsilon +\mathbf{f} &&\qquad\mbox{in }\Ome_T,\\
\Div \bu^\varepsilon &=0  &&\qquad\mbox{in }\Ome_T,\\
\frac{\p \varphi^\varepsilon}{\p t} - \Delta  \mu^\varepsilon +\bu^\varepsilon\cdot \nabla \varphi^\varepsilon &=0 
&&\qquad\mbox{in } \Ome_T,\\
\mu^\varepsilon &=- \eps \Delta \varphi^\varepsilon +\frac{1}{\eps} f(\varphi^\varepsilon)   &&\qquad\mbox{in } \Ome_T.
\end{alignat}
Note that $\varphi^\varepsilon$ is used to denote the phase function in the model. 

Finally, the phase field formulation for the two-phase Hele-Shaw flow model described in {\bf Example 9} of 
Section \ref{sec-1.2.2} is given by the 
following coupled Darcy and Cahn-Hilliard system \cite{CHW17,Feng_Wise12, DFW14, wang_zhang13, han_wang18}: 
\begin{alignat}{2} \label{Darcy_CH} 
 \bu^\varepsilon &= - \nabla p^\varepsilon - \gamma \varphi^\varepsilon \nabla\mu^\varepsilon   &&\qquad\mbox{in }\Ome_T,\\
\Div \bu^\varepsilon &=0  &&\qquad\mbox{in }\Ome_T,\\
\frac{\p \varphi^\varepsilon}{\p t} - \Delta  \mu^\varepsilon +\bu^\varepsilon\cdot \nabla \varphi^\varepsilon &=0 
&&\qquad\mbox{in } \Ome_T,\\
\mu^\varepsilon &=- \eps \Delta \varphi^\varepsilon +\frac{1}{\eps} f(\varphi^\varepsilon)   &&\qquad\mbox{in } \Ome_T.
\end{alignat}
Again, $\varphi^\varepsilon$ is used to denote the phase function.

\begin{remark}
In \cite{huisken_ilmanen01} Huisken and Ilmanen proposed a sub-level set formulation for 
the inverse mean curvature flow (IMCF), as expected, this formulation does not have the form of
the general level set equation \eqref{level_set_1}, which would give the following unusual PDE: 
\begin{equation}\label{IMCF_LS}
\mbox{div}\left(\frac{\nabla u}{|\nabla u|}\right)\cdot \frac{\p u}{\p t} =|\nabla u|. 
\end{equation}

To the best of our knowledge, no phase field formulation for the IMCF was proposed in the literature,
however, for a phase function $u^\eps$ there holds the approximate mean curvature formula
\cite{bellettini13, BMO14, du2005retrieving}:
\begin{equation}\label{mc_formula}
H^\eps = \frac{1}{2} \Bigl( \Delta u^\eps - \frac{1}{\eps^2} f(u^\eps ) \Bigr).
\end{equation} 
Using this formula, the chain rule and the geometric law  for the IMCF, we then get 
the following phase field formulation for the IMCF: 
\begin{equation}\label{PF_IMCF}
\frac{1}{2} \Bigl( \Delta u^\eps - \frac{1}{\eps^2 } f(u^\eps ) \Bigr)\cdot \frac{\p u^\eps}{\p t}
=|\nabla u^\eps|,
\end{equation}
which is also an unusual PDE. 
\end{remark} 

Finally, we like to point out that similar to the formulation of the Allen-Cahn equation as an $L^2$-gradient 
flow discussed earlier 
in  \eqref{phase_field_5} and \eqref{phase_field_6}, many of the phase field models
presented here can also be obtained by the energetic approach. For example 
it is well-known 
\cite{fife88, alikakos94, chen96} that the Cahn-Hilliard equations given in {\bf Example 6} 
is the $H^{-1}$-gradient flow of the Cahn-Hilliard free energy 
$\widehat{\mathcal{J}}:=\varepsilon\mathcal{J}$. That is, 
\begin{align}\label{phase_field_5ch}
\frac{\p u}{\p t} = -\widehat{\mathcal{J}}^\prime (u) \qquad\mbox{in } H^{-1}(\Omega),
\end{align}
where the $H^{-1}$ inner product is defined by
$$
(u,v)_{H^{-1}}=\bigl((-\Delta)^{-1}u, v \bigr)_{L^2},
$$
where $(\cdot, \cdot)_{L^2}$ denotes the standard $L^2$-inner product, and
$(-\Delta)^{-1}$ denotes the inverse of the negative Laplacian (subject to appropriate 
boundary and normalization conditions).

The gradient flow structure implies in particular the dynamic energy law
\begin{align}\label{phase_field_energy_law}
\frac{\p \widehat{\mathcal{J}}(u)}{\p t} = - \left\|\frac{\p u}{\p t} \right\|^2, 
\end{align}
where the $\|\cdot\|$ denotes the norm induced by the inner-product that defines the gradient flow.

In the same fashion, the degenerate Cahn-Hilliard equation  given in  \eqref{phase_field_energy_law}
 can also be interpreted as a weighted $H^{-1}$-gradient 
flow of the free energy $\widehat{\mathcal{J}}$ with the energy law
\begin{align}\label{phase_field_energy_law_dege}
\frac{\p \widehat{\mathcal{J}}(u)}{\p t} = -  \|b(u) \nabla \mu \|_{L^2}^2,
\end{align}
where $\mu:=-\varepsilon \Del u +\frac{1}{\vepsi}f(u)$ is often called the {\em chemical potential}.
On the other hand, equation \eqref{Willmore_PF} in {\bf Example 5}
is the $L^2$-gradient flow of the Willmore energy functional \eqref{Willmore_energy}.
The energy law for the coupled 
Navier-Stokes and Cahn-Hilliard system in {\bf Example 8} can also be formulated using the idea of 
flow map for the total energy that 
is the sum of the phase field interfacial energy and the fluid kinetic energy.

%%%%%%%
\subsection{Relationships Between Phase Field and Other Formulations}
Among five mathematical formulations/methodologies for moving interface problems, the first four, namely, 
the parametric formulation, the (Brakke's) varifold formulation, the (De Georgi's) barrier formulation 
and the level set formulation, all intend to represent and to seek the sharp interface exactly. For 
this reason we call them {\em sharp interface formulations/methodologies}.  In contrast,  the 
phase field formulation and methodology does not represent the interface exactly but only seeks  an approximate 
interface, it introduces a width to the interface so the interface is diffused.  For this reason, it is often
called a {\em diffuse interface} formulation/methodology. So the phase field methodology is 
fundamentally different from other methodologies. 
On the other hand, if we use how the interface is represented and sought as a yardstick, we
also can divide the five formulations/methodologies differently as {\em direct} and {\em indirect}
formulations/methodologies. 
 The parametric, the varifold and the barrier's formulations/methodologies 
belong to the direct approach camp because they represent and seek the interface directly, 
while the level set and the phase field methodologies belong to the indirect approach camp
because they both represent and seek the interface indirectly as a level set of an auxiliary function
(i.e., the level set function and phase field function/variable). 
Comparative studies of phase field methods and other methods are of much interests, see for example, a comparison of phase field and parametric
front tracking reported in \cite{bgn3}. 

Indeed, the level set method and the phase field method are intimately connected. 
Their connection is revealed in the limiting process as the parameter $\eps\to 0^+$
in the phase field formulations/equations.  It can be shown (rigorously or formally) that the diffuse 
interface converges (in some sense) to the sharp interface defined by the level set formulation 
\cite{deMottoni_schatzman95, evans_soner_souganidis92, chen92, Pego89, alikakos94, chen96, abels_lengeler14, abels15, abels_liu18}. 
It should be noted that the convergence may hold beyond onset of singularities. This is quite important because 
it says that theoretically the level set formulation provides a nice mathematical framework and a good solution 
notion for the underlying moving interface problem which the phase field method aims to approximate and to solve, 
although numerically both methods provide different formalisms and platforms for developing efficient 
numerical methods for moving interface problems. 

On the other hand, the level set method and the phase field method are fundamentally different methodologies
because they are based on different principles/ideas (i.e., sharp interface vs diffuse interface). 
The difference also reflect
on their respective PDE models. The level set PDEs often have stronger nonlinearities than their phase field 
counterparts. For example, the level set PDE for the mean curvature flow, see \eqref{level_set_2}, is a 
non-divergence form quasilinear PDE whose solutions are often defined in the viscosity sense \cite{chen_giga_goto91,
	evans_spruck91}, 
while the phase field PDE for the mean curvature flow is a semilinear PDE with constant coefficient 
whose solutions are defined simply by using integration by parts.  
Compared to the phase field method, the level set method has advantages of being simple to 
obtain the level set equation and to provide sharp representation of the interface. 
On the other hand, it often does not satisfy mass conservation, which may lead to the nonphysical 
motions of the interface, and strong nonlinearity of the level set PDE makes it harder to 
solve numerically. Since the phase field function usually has physical meaning, which makes 
the phase field method more convenient to handle physically or biologically charged 
interface motion. However, although the phase field equation is easier to discretize, 
to resolve the thin diffuse interface one is forced to use adaptive mesh techniques 
\cite{PDG99, provatas2011phase,KNS04, Feng_Wu05, Feng_Wu08} 
because the computation becomes intractable, especially in three-dimensional cases if 
uniform meshes are used.
Both the level set method and the phase field method share an important advantage over the 
direct methods, that is, they can handle with ease singularities (or topological changes)
of the interfaces.  

We note that the phase field formulation also relates to other sharp interface 
formulations.  In \cite{ilmanen93} Ilmanen proved the convergence of the solution of a phase field 
equation to the Brakke's varifold solution for the mean curvature flow, and in \cite{bellettini13,bellettini_novaga97} 
Bellettini etc. established the convergence of the solution of another phase field equation to the 
De Georgi's barrier solution for the Willmore flow. 

%%%%
\subsection{Phase Function Representations of Geometric Quantities} \label{sec-1.2.7} 
To obtain a phase field formulation for a given geometric law, it is very important to be 
able to represent well-known geometric quantities such as the mean curvature and
Gauss curvature in terms of the phase function, this then requires 
to develop so-called the phase field calculus (cf.  \cite{du2011phase}).

First, it is well-known that the volume formula of a geometric domain in terms of a phase function $u^{\eps}$ 
can be derived from the integral of the phase field $u^{\eps}$, or more generally, the integral of
$g^\eps(u^\eps) $
where $g^\eps$ could take on an approximation of indicator function such that $g^\eps(u^\eps) $ is
approximately one inside the geometric domain and zero outside (as determined by the values of $u^\eps$).

Second, the interface area formula can be expressed, as discussed before by the Cahn-Hilliard energy functional 
given in \eqref{CH_functional}.

Other frequently used geometric quantities include the interface normal, expressed by $\nabla u^\eps/|\nabla u^\eps|$, and
as already mentioned above, the formula for the mean curvature of the interface represented
by the zero-level set of the phase field function $u^\eps$ is given in \eqref{mc_formula} \cite{du2004phase,du2006simulating}.
In addition, the phase field version
of Gauss curvature can also be computed \cite{bm2010,bellettini13,BMO14,du2005retrieving}. A form similar to \eqref{mc_formula}  is given in  \cite{bm2010} as
$$
%\begin{equation}\label{KC_formula}
K^\eps = \frac{1}{2\eps} \left[ 
\Bigl(  \eps
\Delta u^\eps - \frac{1}{\eps} f(u^\eps ) \Bigr)^2 -
\left|  \eps
\nabla^2 u^\eps - \frac{1}{\eps} f(u^\eps ) 
\frac{\nabla u}{|\nabla u|}
\otimes
\frac{\nabla u}{|\nabla u|} \right|^2  \right],
$$
%\end{equation}
where $|\cdot|$ denotes the standard Euclidean norm.
Consequently, the phase field Euler-Poincar\'e index 
can be computed by integrating the Gauss curvature to retrieve
topological information of the implicitly defined zero level set
surfae of the phase field function \cite{du2005retrieving}. 
A simple form in 2D can be found in  \cite{du2005retrieving} and is given by 
\begin{equation}\label{Euler_formula_2}
\chi^\eps = \frac{C_2}{\eps} \int_\Omega
\Bigl( \Delta u^\eps - \frac{1}{\eps^2} f(u^\eps ) \Bigr) 
d\Omega,
\end{equation}
with the 3D version given in \cite{dlrw2007} by 
$$
%\begin{equation}\label{Euler_formula_3}
\chi^\eps = \frac{C_2}{\eps} \int_\Omega
\Bigl( \Delta u^\eps - \frac{1}{\eps^2} f(u^\eps ) \Bigr) p(u^\eps) d\Omega
%\end{equation}
$$
for a suitably chosen function $p$ such as $p(t)=-2n(1-t^2)^{n-1}t$, $n\geq 1$ and normalization 
constants $C_2$ and $C_3$. Other studies concerning topological constraints in phase field models 
can be found in \cite{dlw17}

%%%%%
\subsection{Convergence of the Phase Field Formulation}\label{sec-1.2.8} 
Since the phase field method is based on a diffuse interface idea/approach, a fundamental 
mathematical question is whether the diffuse interface converges to a sharp interface 
when the width of the diffuse interface tends to zero (i.e., $\varepsilon\to 0$).  
It turns out that this is very difficult question to rigorously answer for all 
moving interface problems. It remains an open problem for many phase field models.  
To address this question, one first needs to know or to guess what is the 
limiting sharp interface formulation. 

Many people contributed to the convergence proof for the  Allen-Cahn problem to the MCF 
(cf. \cite{kohn_sternberg89,keller89, bronsard_kohn91, deMottoni_schatzman95,chen92}),  
it was finally proved by Evans, Soner and Souganidis in \cite{evans_soner_souganidis92} that
	\[
	\mbox{dist}_H\bigl(\G_t^\vepsi,\G_t\bigr) = O(\vepsi),
	\qquad\mbox{as } \vepsi\rightarrow 0,
	\]
where $\mbox{dist}_H$ denotes the Hausdorff distance between two sets.

The connection between the Cahn-Hilliard equation and the Hele-Shaw problem was first 
formally established by Pego in \cite{Pego89}. It was later proved by Alikakos, Bates and Chen
in \cite{alikakos94} that if 
$u^\vepsi$ and $\mu^\vepsi$ satisfy the Cahn-Hilliard system
\[
\frac{\p u^\vepsi}{\p t}=\Del \mu^\vepsi,\qquad
\mu^\vepsi= -\vepsi \Del u^\vepsi +\frac{1}{\vepsi} f(u^\vepsi),
\]
then 
\[
\mu^\varepsilon\to w,\qquad \mbox{dist}_H\bigl(\G_t^\vepsi,\G_t\bigr) \longrightarrow 0
\qquad\mbox{as } \vepsi\rightarrow 0
\]
before onset of singularities (or topological changes). Where $w$ denotes the solution of 
the Hele-Shaw problem. Chen also proved in \cite{chen96} using the energy method 
the convergence of the solution of the Cahn-Hilliard problem 
to a (very) weak solution of the Hele-Shaw problem.

For {\em the surface diffusion flow} $V_n=-\alpha \Del_\G H$, the convergence was formally 
proved by Cahn, Elliott and Novick-Cohn in \cite{cahn_elliott_cohn96} if $u^\vepsi$ and 
$\mu^\varepsilon$ satisfy the following {\em degenerate} Cahn-Hilliard equation/system:
\begin{align*}
\frac{\p u^\vepsi}{\p t}
&= \Div\bigl( b(u^\vepsi) \nab \mu^\vepsi \bigr), \\
\mu^\vepsi &= -\Del u^\vepsi +\frac{1}{\vepsi^2} f(u^\vepsi),
\qquad b(z):=1-z^2.
\end{align*}

For {\em the generalized Stefan problem}, in \cite{caginalp_chen98} Caginalp and Chen 
proved that if $\varphi^\vepsi$ and $u^\vepsi$ satisfies the following parabolic 
system (which is often called the classical phase field model)
\begin{align*}
\alpha(\vepsi)\frac{\p \varphi^\vepsi}{\p t}&= \vepsi\Del\varphi^\vepsi
-\frac{1}{\vepsi} f(\varphi^\vepsi) + s(\vepsi) u^\vepsi, \\
c(\vepsi)\frac{\p u^\vepsi}{\p t}&= \Del u^\vepsi-\frac{\p \varphi^\vepsi}{\p t}, 
\end{align*}
then
\[
u^\varepsilon\to \phi, \qquad \mbox{dist}_H\bigl(\G_t^\vepsi,\G_t\bigr) \longrightarrow 0
\qquad\mbox{as } \vepsi\rightarrow 0
\]
before onset of singularities (or topological changes) for various combinations of $\alpha(\varepsilon)$ 
and $s(\varepsilon)$. Where $\phi$ denotes the solution (temperature) of the generalized 
Stefan problem. 

Other works concerning the sharp-interface limit of Cahn-Hilliard equations and its variants can
be found in \cite{ABK2018,DD2012,CWX2014,dd2014,DLL18,DP2013,GN2000,LMS2016,Novi2000}.
Recently a convergence proof of the Navier-Stokes-Cahn-Hilliard phase field model to the sharp
interface model for two-phase fluids was carried out in \cite{abels_lengeler14, ALS17, abels_liu18} 
and for a Cahn-Larch\'e phase field model approximating an elasticity sharp interface problem 
in \cite{abels15}. 

We conclude this section by summing up the main points of the above convergence results
and also make a few comments. 
First, the level set method is a {\em sharp interface} method,
while phase field method is a {\em diffuse interface} method.
The latter converges to the former as $\vepsi\rightarrow 0$.
Second, the level set function of each moving interface problem must satisfy
a master equation, which is the following Hamilton-Jacobi equation: 
%$$\varphi_t+V_n |\nab \varphi|=0$$
	\[
	\varphi_t+V_n |\nab \varphi|=0 
	%\quad\mbox{or}\quad \varphi_t+V \cdot \nab \vphi=0
	\qquad \mbox{on } \G_t.
	\]
However, there is no such a general master equation for phase field models. 
Phase field equations are usually not unique and different for 
different moving interface problems. 
Third, level set formulations often give rise to quasilinear or fully nonlinear PDEs while 
phase field formulations result in semilinear or quasilinear PDEs which are 
of singularly perturbed type and involve a small parameter $\vepsi$. 
Lastly, level set functions are purely mathematical objects and may not have physical meaning
while phase field functions often represent physical quantities such as densities and concentrations.

\section{Time-stepping Schemes for Phase Field Models}\label{sec-3}
Simulating the time evolution of phase field models has many important applications such as phase transition,
 microstructure coarsening and  cell motion. Typically the dynamic process governed by phase field models involve a 
 number of stages that exhibit features on multiple time and spatial scales. The main focuses on effective time 
 discretization have largely been devoted to the construction of methods which can offer both stable long time 
 simulations and adequate resolutions of transient phenomena. Preserving,  on the discrete level, some mathematical 
 features  associated with the continuum phase field dynamics (such as the energy law and uniform pointwise bounds 
 on the phase field variable) has taken a center stage in much of the numerical analysis literature. It is also 
 desirable that these feature preserving schemes can provide high level resolution (which is not limited to just 
 a formal high order truncation error) and can be efficiently implemented (from effective solver to scalable algorithms). 
 While a number of classical time integrators have been studied, various splitting and stabilization techniques have been 
 proposed to address these concerns. We refer to, for example, \cite{barrett_blowey_garcke99,chen1998applications,chen_nochetto_schmidt00,du_nicolaides91,du17jcp,dz09numerical,elliott_french87,elliott_french89,eyre1998,Feng_Prohl03b,Feng_Prohl04a,Feng_Prohl05a, Feng_Prohl04c, Feng_Li15, FLX16a, wu_li18, gomez2018computational,li2016,shen2010,shen2015,shen2018,wise2009,YDZ18} and the references cited therein.

We note in particular that past works have utilized the gradient flow structure of many phase field equations 
which we focus on here. As an illustration,  following \eqref{phase_field_5ch}, we formulate a generic phase 
field equation as $u_t = -\mathcal{J}^\prime (u)$ with $u_t$ being the time derivative and $\mathcal{J}^\prime$ 
representing the variation of the phase field energy with respect to a generic inner product (topology) that 
could be either $L^2$ (for Allen-Cahn dynamics) or $H^{-1}$ (for Cahn-Hilliard dynamics). We use $<\cdot,\cdot>$ 
to represent the associated inner product for convenience.

To present time discretization schemes, we let $\{t_n\}_{n=1}^N$ denote the discrete time steps 
with $t_0$ be the initial time and $t_N=T$ be the terminal time.  We use $\tau$ to denote the discrete time 
step size, generically $\tau$ can be dependent on $t_n$ , that is, $\tau=\tau_n$, which gives nonuniform 
time steps or adaptive time steps. For the sake of simple notation, we drop the explicit dependence on $n$ 
for the discussion in this section.  We use $u^n$ to denote the numerical approximation of $u(t_n)$,  the function 
value of the phase field variable $u$ at time $t_n$.

\subsection{Classical Schemes}

By viewing phase field equations as abstract dynamic systems, many classical time discretization schemes can be used.  With the gradient flow formulation, one may offer various interpretations to these classical schemes.

{\bf Example 1:}   ({\em Fully-explicit/forward Euler Scheme}) 
Standard fully explicit Euler schemes for the gradient system $u_t=-\mathcal{J}^\prime (u)$  can be written as
\begin{align}\label{td-1}
u^{n+1}= u^n  -\tau \mathcal{J}^\prime (u^{n}).
\end{align}
It can be viewed as a steepest descent iteration of the phase field energy, in the corresponding topology, 
with a step size $\tau$. While popular among domain scientists due to the simplicity in the numerical implementation, 
the limitation on the time step size, due to stability considerations, severely affects its performance for long time integration.  

{\bf Example 2:}   ({\em Fully-implicit/backward Euler Scheme})  To deal with the stiff nature of phase field models, 
implicit schemes that often better stability properties are better choices for long time integration than the conventional 
explicit Euler scheme. 

The standard backward/implicit Euler scheme can be written as
\begin{align}\label{td-2}
u^{n+1}= u^n  -\tau \mathcal{J}^\prime (u^{n+1}),
\end{align}
which may be reformulated in a variational form as
$$\min \Bigl\{   \mathcal{J}(u)+\frac{1}{2\tau} \|u-u^n\|^2 \Bigr\},$$
where $\|\cdot\|$ is the norm corresponding to the same inner product $<\cdot, \cdot>$ that defines the variational derivative of $J$.

Thus, by interpreting $1/\tau$ as a penalty constant, one can view the backward/implicit Euler scheme as to look for the minimum of $J$ near $u^n$. A consequence is that
$$
  \mathcal{J}(u^{n+1}) +\frac{1}{2\tau} \|u^{n+1}-u^n\|^2
=
\min \Bigl\{ \mathcal{J}(u)+\frac{1}{2\tau} \|u-u^n\|^2 \Bigr\} \leq  \mathcal{J}(u^n).
$$
In particular, this leads to the decrease of energy for any step size.  Nevertheless, to compute $u^{n+1}$, a non-convex problem needs to be solved and condition of $\tau$ may need to be imposed to ensure the uniqueness of the minimizer.

{\bf Example 3:}   ({\em  Crank-Nicolson and its variant})
For second order Crank-Nicolson type schemes, we may consider
$$
u^{n+1}= u^n  -\tau \mathcal{J}^\prime \left(\frac{u^{n+1}+u^n}{2} \right),
$$
which may be  obtained via the  extrapolation $2 u^{n+1}_*- u^n$ where $u^{n+1}_*$ is backward Euler solution. Other variants include
$$
u^{n+1}= u^n  -\frac{\tau}{2} \mathcal{J}^\prime (u^{n+1}) -\frac{\tau}{2} 
 \mathcal{J}^\prime (u^n),
$$
which solves 
$$
\min \Bigl\{ \mathcal{J}(u)+\frac{1}{\tau}\|u-u^n\|^2 + \bigl(u-u^n,  \mathcal{J}^\prime (u^n) \bigr) \Bigr\},
$$
instead, and a modified scheme
 \begin{align}\label{td-3}
u^{n+1}&= u^n  -\tau \phi(u^n,u^{n+1}), \\
<\phi(v,w), v-w> &=  \mathcal{J} (v) -
 \mathcal{J} (w) \qquad \forall v, w.
\end{align}
For the latter, we have
$$\frac{1}{\tau} 
\Bigl( \mathcal{J} (u^{n+1}) -  \mathcal{J} (u^n) \Bigr) + \frac{1}{\tau} \| u^{n+1}-u^n\|^2 = 0.
 $$
 While for other schemes one can often only derive an inequality version the energy law, the above
  can be seen as an exact time discrete analogue of the energy law associated with the continuum model.

The use of the discrete variation $\phi(v,w)$ has been considered in other works such 
as \cite{du2011finite,furihata2001stable}. Different variants have also been studied based on Taylor 
expansions \cite{gomez2011provably}.

 Besides classical schemes considered above, there are also works on numerical discretization of phase 
 field models based on traditional linear multistep methods   \cite{akrivis1998} and Runge-Kutta (RK)
 methods \cite{song2017,shin2017unconditionally} which also can be studied for general gradient 
 system  \cite{humphries1994runge}.
 
 \subsection{Convex Splitting and Stablized Schemes}
 It is desirable to have time discretization of the gradient flows satisfying the energy stability which refers to
 $$\mathcal{J}(u^{n+1})\leq \mathcal{J}(u^n).$$
   In particular, it is often viewed that the decay of energy not only preserves  the physical feature 
   but also provides stability to numerical discretization.  For some phase models, preserving the pointwise 
   bound enjoyed by the continuum model is also an important feature of numerical schemes.  While a number 
   of classical schemes mentioned in the above share such properties, there are more variants that offer
    more efficient implementation and better performance.
 
 {\bf Example 1:} {\em Convex splitting}
 
The fully implicit schemes enjoy energy stability but at the expense of solving nonlinear non-convex 
problems in general. An alternative is the convex splitting for gradient flows \cite{elliott1993global,eyre1998}. 
Here, based  on a decomposition $J=J_1+J_2$ where $J_1 $ is convex and $J_2$ is concave, the convex splitting scheme is given by
\begin{align}\label{td-4}
u^{n+1}= u^n  - \tau \mathcal{J}^\prime_1 (u^{n+1})
- \tau \mathcal{J}^\prime_2 (u^{n}),
\end{align}
which also has an equivalent variations formulation
$$\min \Bigl\{ \mathcal{J}_1(u)+\frac{1}{2\tau}\|u-u^n\|^2-  \bigl( u-u_n,  \mathcal{J}^\prime_2 (u^n) \bigr) \Bigr\}. $$
The convexity of  $\mathcal{J}_1$ assures the unique solution for any $\tau>0$. The energy stability also follows,
that is, 
$$\mathcal{J}(u^{n+1})+\frac{1}{2\tau}\|u^{n+1}-u^n\|^2\leq  
\mathcal{J}(u^{n}).
$$

There are a number of higher order extensions, for example, by combining the convex splitting with Runge-Kutta 
methods.

{\bf Example 2:}   {\em Linearly-implicit stabilized schemes}.

To make the computational task simpler for each time step than even that associated with the convex splitting schemes,  
a popular technique for getting  energy stable scheme is to  use a semi-implicit or  coupled implicit-explicit scheme 
that lead to only linear systems of equation.

Among the first order scheme, this generically amounts to a linear splitting of the $ \mathcal{J}^\prime (u)$ into 
the sum  $\mathcal{J}^\prime (u) - Au=N(u)$ and $ Au$ for some linear operator $A$ where the first term is treated 
explicitly while the second is treated implicitly. That is,
\begin{align}\label{td-5}
(I+ \tau A) u^{n+1}= u^n  - \tau N(u^n),
\end{align}
Typically $Au$ is taken to be a scalar  linear combination of the  variation of the quadratic energy 
involving $|\nabla u|$ and $u$ itself. The scalars in the combination are choosing to make $A$ sufficiently 
coercive so that \eqref{td-5} becomes dissipative and energy stability can be assured.

The energy stability and the need for only linear solvers per time step have made the first order stabilized scheme 
 popular in practice. Similar work on second order BDF scheme can be found in 
\cite{YCWW18}. 
For studies on other higher order stabilization, we refer to  \cite{li2017,song2017}.

{\bf Example 3:} { \em Exponential integrator}
       
 Another way of utilizing  the   linear splitting of the $ \mathcal{J}^\prime (u)=Au+N(u)$ is to develop exponential integrators \cite{cox2002,hochbruck2010,kassam2005}, The central idea is to use an equivalent formulation of the phase field equation 
\begin{align}
\label{td_exp1}
u(t+\tau)=e^{\tau A }u(t)+\int_0^\tau  e^{(\tau-s)A} N(u(t+s))\,\d s,\
\end{align}
or
\begin{align}
\label{td_exp2}
e^{- \tau A }u(t+\tau)=u(t)+\int_0^\tau  e^{sA} N(u(t+s))\,\d s.
\end{align}
By applying different quadratures to the integral over time, various discrete integration schemes can be derived. 
 For \eqref{td_exp1}, if the quadratures are based polynomial interpolations of $N(u(t+s))$ using the discrete values 
 of $u$ from previous time steps, we end up with the ETD (Exponential Time Difference) linear multistep schemes. 
 Meanwhile, the ETD-RK type scheme 
 is to approximate $N(u(t+s))$ via polynomial interpolations of $\{N(u(t+s_i))\}$ where $\{s_i\}$ are the quadrature 
 nodes in $(0,\tau)$ and the nodal values are constructed stage-wise through the lower order ETD-RK schemes.  On the 
 other hand, if we consider \eqref{td_exp2},  a polynomial interpolation can be done for the entire integrand 
 $ e^{sA} N(u(t+s))$ instead, which would lead to the IF-RK schemes if a multistage interpolation is used.

As an illustration, the ETD-RK1 scheme is simply
\begin{align}
\label{td_ETD1}
u^{n+1}=e^{\tau A} u^n+   (e^{\tau A}-I) A^{-1} N(u^n),
\end{align}
and the ETD-RK2 scheme is given by
\begin{align}
\label{td_ETD2}
u^{n+1}=\tilde{u}^{n+1} + \tau^{-1}  \Big(e^{-\tau A}-1+ {\tau}A \Big) 
A^{-2} \Big[ N (\tilde{u}^{n+1}) - N(u^n)\Big]\,
\end{align}
where $\tilde{u}^{n+1}$ is obtained from \eqref{td_ETD1}.

For both schemes. stability can be established for various phase field models \cite{du05etd,du05etd,kassam2005,djlq18}. Long time integrations and
large scale simulations have been carried out \cite{jzd2015,zhang16,wang2016efficient}.

We note that if a first order Taylor expansion of the exponential is used, then we recover the first 
order   explicit Euler scheme if  $e^{\tau A} \sim I +\tau A$, or the first order  semi-implicit 
Euler $( I-\tau A) u^{n+1}=  u^n + \tau  N(u^n)$ if  $e^{-\tau A} \sim I-\tau A$.

	\subsection{Schemes Using Lagrangian Multipliers}
	
	Recently an interesting approach was proposed in \cite{guillen_tierra13} for constructing 
	unconditional stable linear schemes, it is based on a Lagrange multiplier approach introduced in 
	\cite{BGG11}. This approach has lately been extensively developed and expanded, two new families of time-stepping 
	schemes have been obtained,  
	namely, {\em invariant energy quadratization} (IEQ)  \cite{yang16,zhao_yang_shen_wang16} 
	and {\em scalar auxiliary variable schemes} (SAV) schemes \cite{shen2018}, which 
	are applicable to a large class of free energies. Below we explain the ideas 
	of each family using the constant mobility Cahn-Hilliard equation \eqref{cahn_hilliard} 
	as an example. 
	
	{\bf Example 1:} { \em Invariant energy quadratization schemes}
	
	Consider the 
	Cahn-Hilliard equation \eqref{cahn_hilliard} and suppose that $F(u)\geq -C_0$ for
	some positive constant $C_0$ (note that $F'(u)=f(u)=u^3-u$ and $C_0=0$ for the Cahn-Hilliard equation). 
	Let $q:=\sqrt{F(u)+C_0}$ and rewrite \eqref{cahn_hilliard} as the system 
	\begin{align*}
	u_t &=\Del \mu, \\
	\mu &=- \varepsilon \Del u +  \frac{1}{\varepsilon} G(u) q, \\
	q_t &=\frac12 G(u)  u_t,  \qquad G(u):= \frac{f(u)}{\sqrt{F(u)+C_0}}.
	\end{align*}
	Based on the above reformulation, the following linear implicit-explicit (IMEX) scheme, which is referred
	as an IEQ scheme, can be easily obtained:
	\begin{subequations}\label{CH_IEQ}
		\begin{align}
		u^{n+1}-u^n &= \tau \Del \mu^{n+1}, \label{CH_IEQa} \\
		\mu^{n+1} &=- \varepsilon \Del u^{n+1} +  \frac{ 1 }{\varepsilon} G(u^n)  q^{n+1}, \label{CH_IEQb}\\
		q^{n+1}-q^n &=\frac12 G(u^n)  \bigl(u^{n+1}-u^n \bigr). \label{CH_IEQc}
		\end{align}
	\end{subequations}
	
	It is easy to check that the above IEQ scheme has $O(\tau)$ truncation error and it is unconditional 
	energy stable for the discrete energy functional
	\[
	E^n:= \frac{\varepsilon}2 \|\nabla u^n\|_{L^2}^2 + \varepsilon \|q^n\|_{L^2}^2. 
	\]
	It needs to be noted that the above discrete energy functional is a modification of the original energy functional
	and both are the same at the PDE level. For its detailed analysis and higher oder extensions as well 
	as applications to other problems, we refer to \cite{yang2017linear, yang19} and the references therein.

	{\bf Example 2:} { \em Scalar auxiliary variable schemes}
	
	Notice that although the above IEQ only solves a linear system (after a spatial discretization 
	is applied) at each time step, its coefficient 
	matrix varies at each time step.  To remedy this drawback, 
	%(and a couple of other disadvantages of IEQ schemes), 
	the SAV approach starts with the induction of an auxiliary scalar function  
	\[
	r:=\sqrt{E_1(u) +C_1}, \qquad\mbox{where}\quad E_1(u):=\int_\Omega F(u)\, dx,
	\]
	assuming that $E_1(u)\geq -C_1$. Then the original Cahn-Hilliard equation can be 
	rewritten as the following system:
	\begin{align*}
	u_t &=\Del \mu, \\
	\mu &=- \varepsilon \Del u +  \frac{1}{\varepsilon}  H(u) r,\\
	r_t &=\frac{1}{2} \bigl( H(u), u_t \big)  ,\qquad  H(u):=\frac{f(u)}{\sqrt{E_1(u)+C_1} }.
	\end{align*}
	Based on the above reformulation, the following linear IMEX scheme, which is referred 
	as a SAV scheme, can be easily obtained:
	\begin{subequations}\label{CH_SAV}
		\begin{align}
		u^{n+1} - u^n &= \tau \Del \mu^{n+1}, \label{CH_SAVa}\\
		\mu^{n+1} &= -\varepsilon \Del u^{n+1} + \frac{1}{\varepsilon} H(u^n) r^{n+1}, \label{CH_SAVb}\\
		r^{n+1} -r^n & =  \frac{1}{ 2  }  \bigr( H(u^n), u^{n+1}-u^n \bigr), \label{CH_SAVc}
		\end{align}
	\end{subequations}
	
{ It should be noted that constant coefficient equations only need to be solved 
		at each time step because auxiliary variable $r^{n+1}$ can be eliminated}. 
	    It is easy to check that the above SAV scheme has $O(\tau)$ truncation 
	error and it is unconditional energy stable for the discrete energy functional
	\[
	G^n:= \frac{\varepsilon}2 \|\nabla u^n\|_{L^2}^2 + \varepsilon \bigl(r^n\bigr)^2. 
	\]
	It also needs to be pointed out that the above discrete energy functional is a modification of the original 
	energy functional and both are the same at the PDE level. For its detailed analysis and higher oder 
	extensions as well as applications to other problems, we refer to \cite{cheng_shen18, shen_xu18, sxy2018}
	and the references therein.

\subsection{Further Considerations}
To design  effective time integrator, 
there are other important aspects besides the
construction of the time discretization.
For  example, spatial approximations require 
particular attention due to the appearance of
high order spatial derivatives in many phase field models
such as the Cahn-Hilliard equation, see discussions in Section \ref{sec-4}.
A related matter is the design of efficient linear and nonlinear solvers to find 
the discrete solutions at each time step (from
Newton's methods, Quasi-Newton, optimization based methods,
preconditioner techniques, multilevel methods, and domain decomposition
and scalable algorithms). 
While some time discretizations lead to difficult nonlinear systems, 
others (such as ETD schemes) avoids the linear solver completely.
The latter, when combined with domain decomposition techniques, can offer
high scalability  as demonstrated recently \cite{ju2015fast,zhang16}.

Concerning long time integration, there is much potential to get  more effective simulations 
via variable time steps. This is discussion in the context of adaptive methods in Section \ref{sec-6}.
Another aspect is the study of time integrators for coupling phase field and
other physical  models such as Cahn-Hilliard-Navier-Stokes equations and
Cahn-Hilliard-Microelasticity models. This will be discussed later in Section \ref{sec-7}.
We also note that the time integration  of stochastic phase field models such as stochastic 
Allen-Cahn and Cahn-Hilliard equations are also topics studied in the litertature, we refer the reader 
to \cite{OWW14, weber_roger13, FKLL18, KLL15, KLM14} and the references therein for detailed 
discussions.

Related to the subject, there are algorithmic issues in search of transition states and critical nuclei 
of the phase field energy, which can help shed light on processes of phase transformation. Instead of 
gradient dynamics that lead to equilibria, gentlest ascend dynamics and shrinking dimer dynamics have 
been developed for converging to saddle points \cite{zhang2012shrinking,zhou20103d}. Various optimization algorithms 
and time stepping schemes have also been developed \cite{zhang2016recent}.

\section{Spatial Discretization Methods for Phase Field Models}\label{sec-4}

Spatial discretization of phase field models have been widely studied, ranging from conventional
finite difference, finite element and spectral methods,  to specialized techniques 
like radial basis functions and isogeometric analysis,  In this section, we present the basic formulations 
of a number of spatial discretization schemes.

\subsection{Spatial Finite Difference  Discretization}\label{sec-4.1}
For numerical solution of PDE models, finite difference methods are often the ones studied first before other discretization techniques. Along with the discussion of phase field modeling and simulations of free boundary problems in the 1980s, there have been works on their spatial approximations by finite difference methods \cite{caginalp1987numerical,lin1988numerical}. Much of the difference approximations are standard, for example, second order center differences are used to approximate the second order spatial derivatives. Higher order spatial derivatives are then approximated by the compositions of differences.

For example, for the Allen-Cahn equation \eqref{phase_field_5a} 
and the Cahn-Hilliard equation \eqref{cahn_hilliard} with a constant mobility, 
the semi-discrete in space approximation can be given respectively by 
$$
\frac{\p u^h_j}{\p t} = (\Delta_h u^h)_j - \frac{1}{\eps^2} f(u^h_j).
$$
and
$$
\frac{\p u^h_j}{\p t} +\Del_h \Bigl( \vepsi \Del_h u_j^h  -\frac{1}{\vepsi} f(u^h_j) \Bigr) =0,
$$
Here $u^h_j$ denotes the numerical solution at spatial grid point (labeled by the
subscript $j$) with a typical uniform grid spacing $h$. $\Del_h$ denotes the
standard discrete Laplacian based on the second order center difference.
Coupled with time integration schemes, one can get fully discrete approximation, see for example the work in \cite{du_nicolaides91} on the Cahn-Hilliard equation that can preserve the energy law at the discrete level. Further numerical analysis work on the difference approximations can also be found in \cite{sun1995second,furihata2001stable}. Analysis of difference approximations in space under periodic conditions coupled with some BDF energy stable time discretization 
was presented in \cite{cheng2018energy}. For the Allen-Cahn equation,  the analysis of difference approximations  can be
found in \cite{chen1998convergence}.

Starting from the earlier attempts, finite difference approximations have been adopted in many subsequent computational investigations of spinodal decomposition, dendritic growth, solidification and so on, see for example, \cite{rogers1988numerical}. In particular, there are discussions, largely based on numerical experiments, on the choices of time steps and spatial grid size with respect to the diffuse interfacial width parameter $\eps$, see for example, \cite{wheeler1993computation}. In the application domain, difference approximations remain popular for many practical applications of phase field methods in recent years as exemplified by the Gordon Bell prize winning work \cite{GB2011}.

On the algorithm development and numerical analysis side, more recent focus has been on effective linear and nonlinear  solvers, including the development of various adaptive and multigrid techniques, see \cite{kim2004conservative,rosam2008adaptive,wise2007solving} and discussion in later sections. High order compact difference spatial approximations have also been studied \cite{li2016compact} which, when coupled with high order time integrators, have led to record-breaking extreme scale 3D simulations of coarsening dynamics based on the Cahn-Hilliard model.

\subsection{Spatial Galerkin Discretizations}\label{sec-4.2}
%The finite difference method is based on the idea of approximating derivatives in 
%a PDE by finite difference quotients. Such an approximation is accurate only when the 
%underlying solution function is sufficiently smooth. This can be easily seen from 
%error estimates for finite element solutions which involve high order classical 
%derivatives of the PDE solution. 
%
Unlike the finite difference method that is based on the idea of approximating derivatives in 
a PDE by discrete finite difference operators, another larger class of numerical methods, called 
the Galerkin method, for PDEs are based on a very different idea, namely, to employ 
a variational principle (or a weak formulation) and to approximate an infinite-dimensional 
Banach (or Hilbert) space by a sequence of finite-dimensional spaces which may or may not 
be subspaces of the infinite-dimensional Banach space. 

Again, we use the Allen-Cahn equation \eqref{phase_field_5a} and the Cahn-Hilliard 
equation \eqref{cahn_hilliard} as two examples to
demonstrate the formulation of their Galerkin approximations and the widely used 
four types of Galerkin methods, that is, traditional finite element methods, spectral methods,
discontinuous Galerkin methods, and isogeometric analysis.  

A simple integration by parts immediately leads to the following weak formulation for
the Allen-Cahn equation \eqref{phase_field_5a} complemented by the Neumann boundary 
condition: Find $u:(0,T)\to V:=H^1(\Omega)$ such that 

\begin{align}\label{AC_weakform}
\bigl(u_t, v\bigr) + a\big( u, v\bigr) + \frac{1}{\varepsilon^2} \bigl(f(u), v\big)
=0\qquad \forall v\in V,
\end{align}  
where and throughout this section $(\cdot,\cdot)=(\cdot,\cdot)_\Omega$ denotes 
the $L^2(\Omega)$-inner product and 
\[
a\big(u, v\bigr) := \big(\nabla u, \nabla v\bigr). 
\]

For the Cahn-Hilliard equation \eqref{cahn_hilliard} with no-flux boundary conditions,
its weak formulation is defined as finding $u:(0,T)\to \widetilde{V}=H^2(\Omega)$ such that 

\begin{align}\label{CH_weakform}
\bigl(u_t, v\bigr) + \varepsilon \widetilde{a}\big( u, v\bigr) 
+\frac{1}{\varepsilon}  c(u, v) =0\qquad \forall v\in  \widetilde{V},
\end{align}  
where 
\[
\widetilde{a} \big( w, v\bigr):= \big(\Delta w, \Delta v\bigr), \qquad
c\bigl(w, v\big):= a\bigl(f(w), v\big).
\]

Let $N (>>1)$ be a positive integer and $V_N$ denote an $N$-dimensional approximation of 
$V$.  In addition, let $a_N(\cdot,\cdot)$, which is defined on $V\cup V_N\times V\cup V_N $, 
denote an approximation 
of $a_N(\cdot,\cdot)$, then the (semi-discrete in space) Galerkin method 
for \eqref{AC_weakform} is defined as seeking $u_N: (0,T)\to V_N$ such that 
\begin{align}\label{AC_GalerkinApprox}
\bigl((u_N)_t, v_N\bigr) + a_N\big( u_N, v_N\bigr) + \frac{1}{\varepsilon^2} \bigl(f(u_N), v_N\big)
=0\qquad \forall v_N\in V_N.
\end{align}  
Similarly, 
the (semi-discrete in space) Galerkin method 
for \eqref{CH_weakform} is defined as seeking $u_N: (0,T)\to \widetilde{V}_N$ such that 
\begin{align}\label{CH_GalerkinApprox}
\bigl((u_N)_t, v_N\bigr) + \varepsilon \widetilde{a}_N\big( u_N, v_N\bigr) + \frac{1}{\varepsilon} c_N\bigl(u_N, v_N\big)
=0\qquad \forall v_N\in \widetilde{V}_N.
\end{align}  

To define specific Galerkin methods, one needs to choose some $N$-dimensional spaces $V_N$ 
and $\widetilde{V}_N$, specific bilinear forms $a_N(\cdot.\cdot)$, $\widetilde{a}_N(\cdot.\cdot)$.
and a nonlinear form $c_N(\cdot,\cdot)$. It turns out that the choice of 
the latter often depends on the choice of the former and different combinations of $V_N$ 
(resp. $\widetilde{V}_N$) and $a_N(\cdot,\cdot)$ (resp. $\widetilde{V}_N(\cdot,\cdot)$) lead to different 
Galerkin methods for problem \eqref{AC_weakform} (resp. \eqref{CH_weakform}). 
Below we briefly discuss four types of Galerkin methods which are widely used in 
the literature to approximate phase field models.

\subsubsection{Galerkin Discretization via Finite Element Methods}\label{sec-4.2.1}

To define the finite element spatial discretization,  let $\mathcal{T}_h$ denote a quasi-uniform 
(conforming) triangular or rectangular mesh for the physical domain $\Omega$ and $h (>0)$ denote 
the mesh size. The finite element method can be viewed as a Galerkin method in which the finite
dimensional space $V_N$ (resp. $\widetilde{V}_N$), denoted by $V_h$ (resp. $\widetilde{V}_h$), 
consists of piecewise polynomial functions 
(of a fixed degree) over $\mathcal{T}_h$. Here one may think that $h=O(N^{-\frac{1}{d}})$. 
There are two scenarios to be considered separately. First, 
$V_h \subset V$ (resp. $\widetilde{V}_h\subset \widetilde{V}$). In this case, we define $a_N(\cdot, \cdot):=a(\cdot,\cdot)$ (resp. $\widetilde{a}_N(\cdot, \cdot):=\widetilde{a}(\cdot,\cdot)$)
and the method is called 
a {\em conforming finite element method}.  Second, $V_h\not\subset V$ (resp. 
$\widetilde{V}_h\not\subset \widetilde{V}$ ), the resulting method
is called a {\em nonconforming finite element method}. In this case, 
$a_N(\cdot, \cdot)=a(\cdot,\cdot)$ 
(resp. $\widetilde{a}_N(\cdot, \cdot):=\widetilde{a}(\cdot,\cdot)$) 
does not work anymore because it may not even be defined 
on $V_h\times V_h$ (resp. $\widetilde{V}_h\times \widetilde{V}_h$).  The remedy for coping with this difficulty 
depends on the structure of $V_h$ (resp. $\widetilde{V}_h$), however, the choice 
$a_N(\cdot,\cdot)= a_h(\cdot,\cdot)$ (resp. 
$\widetilde{a}_N(\cdot, \cdot)=\widetilde{a}_h(\cdot,\cdot)$
and $c_N(\cdot,\cdot)=c_h(\cdot,\cdot)$) 
works most of the time, where 
\begin{align} \label{AC_FEM}
a_h(w, v)&:=\sum_{K\in\mathcal{T}_h} a_K(w,v) \qquad\forall v,w \in V_h\cup V, \\
\widetilde{a}_h(w,v)&:=\sum_{K\in\mathcal{T}_h} \widetilde{a}_K(w,v) \qquad \forall v,w 
\in \widetilde{V}_h \cup \widetilde{V}, \\
c_h(w,v)&:= \sum_{K\in\mathcal{T}_h} c_K(w,v) \qquad \forall v,w \in \widetilde{V}\cup \widetilde{V}, 
\end{align}
and $a_K(\cdot,\cdot):=a(\cdot,\cdot)|_K$, $\widetilde{a}_K(\cdot,\cdot):= \widetilde{a}(\cdot,\cdot)|_K$,
and $c_K(\cdot,\cdot):= c(\cdot,\cdot)|_K$, 
they are respectively the restrictions of $a(\cdot,\cdot)$,
$\widetilde{a}(\cdot,\cdot)$, and $c(\cdot,\cdot)$ on the element $K\in\mathcal{T}_h$. 
It should be noted that for the Cahn-Hilliard equation, the construction of both conforming 
and nonconforming finite element spaces $\widetilde{V}_h$ has been a nontrivial problem, 
especially in the $3$D case (\cite{Ciarlet78, hu_zhang17}). 

Conforming and nonconforming finite element methods for the Allen-Cahn equation were studied 
in \cite{chen1998convergence, bates2009numerical, Feng_Prohl03b}.
Conforming and nonconforming finite element methods for the Cahn-Hilliard equation were 
studied in \cite{barrett_blowey99, barrett_blowey_garcke99, elliott_french87, elliott_french89, 
	du_nicolaides91, Feng_Prohl04a, Feng_Prohl05a, du2011finite, wu_li18, 
zhang2010nonconforming}.

\subsubsection{Galerkin Discretization via Spectral Methods}\label{sec-4.2.2}
The spectral element method can be viewed as a Galerkin method in which the finite dimensional 
space $V_N$ (resp.  $\widetilde{V}_N$) is chosen as the space of (global algebraic or trigonometric) polynomials whose degree does not exceed $m (>>1)$ in each coordinate direction. Hence,
$N=(m+1)^d$. In practice, $V_N$ (resp. $\widetilde{V}_N$) is often written as a tensor 
product space and is expanded by 
an orthogonal polynomial basis such as Legendre polynomials and Chebyshev polynomials.  
Since functions in $V_N$ (resp.  $\widetilde{V}_N$) are smooth, there always hold $V_N\subset V$ 
and $\widetilde{V}_N \subset \widetilde{V}$. Hence, all spectral methods are confirming methods and 
one naturally chooses $a_N(\cdot, \cdot)=a(\cdot,\cdot)$ and
$\widetilde{a}_N(\cdot, \cdot):=\widetilde{a}(\cdot,\cdot)$ 
%and $c_N(\cdot,\cdot)=c(\cdot,\cdot)$
 in the general Galerkin framework. 

For problems defined on regular computational domains but not with periodic boundary conditions, other
spectral implementations are possible such as Chebyshev spectral methods for Dirichlet boundary conditions,
see discussions in the following books and recent surveys \cite{boyd2001chebyshev,canuto2007spectral,gottlieb1983numerical,guo1998spectral,luo2016parameter,shen2011spectral,bernardi1997spectral}.

Optimal error estimates can be found
in \cite{shen2010} for the above approximation and a number of other methods involving spectral approximations in space and semi-implicit or fully-implicit in time schemes.
These estimates  show that for a fixed $\varepsilon$, the convergence rate is exponential in terms of
the number of the Fourier modes or the dimension of the associated  finite dimensional spaces.

Spectral methods can be extended to more complex phase field models such as those involving
elastic interactions via microelasticity theory \cite{hu2001phase}.
Adaptive spectral approximations for phase field models have been implemented in the context of moving mesh methods \cite{feng2006spectral}.

\subsubsection{Galerkin Discretization via DG Methods}\label{sec-4.2.3}
Although conceptually discontinuous Galerkin (DG) methods can be viewed as nonconforming 
finite element methods, it is often presented as a distinct class of Galerkin methods because they 
use {\em totally discontinuous} piecewise polynomials, while classical nonconforming finite element
methods use piecewise polynomial functions which are not totally discontinuous across element 
edges/faces. So all DG methods (for second and higher order PDEs) are nonconforming Galerkin methods.   
In addition, as expected, the convenience of using totally discontinuous piecewise polynomial 
functions should add some complications for designing bilinear forms $a_N(\cdot, \cdot)$ and 
$\widetilde{a}_N(\cdot, \cdot)$ as well as the nonlinear form $c_N(\cdot,\cdot)$. 

Again, let $\mathcal{T}_h$ denote a (conforming or nonconforming) triangular 
or rectangular mesh for the physical domain $\Omega$ and $h (>0)$ denote
the mesh size. Also let $\mathcal{E}^I$ and $\mathcal{E}^B$ respectively denote the interior and 
boundary edge/face sets and define $\mathcal{E}:=\mathcal{E}^I \cup \mathcal{E}^B$, 
The finite dimensional space $V_N$ (resp. $\widetilde{V}_N$), 
denoted by $V^{\mbox{\tiny DG}}_h$ (resp. $\widetilde{V}^{\mbox{\tiny DG}}_h=V^{\mbox{\tiny DG}}_h$) 
and called it the DG space, consists of 
piecewise polynomial functions (of a fixed degree) over $\mathcal{T}_h$ which are 
totally discontinuous across element edges/faces. Here one may think that $h=O(N^{-\frac{1}{d}})$.
The choices for $a_N(\cdot, \cdot)$ and $\widetilde{a}_N(\cdot, \cdot)$ are 
more complicated, the following are two widely used DG discretizations of $a(\cdot, \cdot)$ and $\widetilde{a}(\cdot,\cdot)$ \cite{Feng_Li15,Feng_Karakashian07}:
\begin{align*}
a_h(w,v)&:= \sum_{K\in\mathcal{T}_h}   a_K(w,v)  
           - \sum_{e\in \mathcal{E}^I} \Bigl( \Langle \{ \partial_n w\}, [v] \Rangle_e
           + \lambda \bigl\langle [w], \{ \partial_n v\} \bigr\rangle_e  \\
        &\hskip 1.5in  + \frac{\sigma_e}{h_e}\Langle [w],[v] \Rangle_e \Bigr), \\
\widetilde{a}(w,v) &:= \sum_{K\in\mathcal{T}_h}   \widetilde{a}_K(w,v)  
          +\sum_{e\in\mathcal{E}^I} \Big( \Langle \{\partial_n \Delta w\}, [v] \Rangle_e
          + \lambda \Langle \{ \partial_n \Delta v\}, [w] \Rangle_e \\
         &\qquad + \gamma_e h_e^{-3} \Langle [w],[v] \Rangle_e \Big)
          - \sum_{e\in\mathcal{E}} \Big( \Langle \{\Delta w\}, [\partial_n v] \Rangle_e
          + \lambda \Langle \{\Delta v\}, [\partial_n w] \Rangle_e 
\end{align*}
\begin{align*}
    - \beta_e h_e^{-1} \Langle [\partial_n w], [\partial_n v] \Rangle_e \Big)
\end{align*} 
for all $v,w \in V^{\mbox{\tiny DG}}_h$. Where $\langle\cdot,\cdot\rangle_e$ stands 
for the $L^2$-inner product on $e$,
$\partial_n v$ denotes the normal derivative (on $e$), $[\cdot]$ and $\{\cdot\}$ denote
respectively the jump and average operators (on $e$). Moreover, $\sigma_e, \beta_e$ and $\gamma_e$ 
are under-determined positive constants, and $\lambda=\{-1,0,1\}$ with $\lambda=1$ giving  
two symmetric DG bilinear forms. 

The choices for $c_N(\cdot,\cdot)$ are various. The simplest one is
\[
c_h(w,v):= \sum_{K\in\mathcal{T}_h} c_K(w,v)  \qquad \forall v,w\in V^{\mbox{\tiny DG}}_h. 
\]
Another more involved choice \cite{Feng_Karakashian07} is 
\begin{align*}
c_h(w,v) := -\sum_{K\in\mathcal{T}_h} \bigl( f(u), \Delta v \bigr)_K 
&+\sum_{e\in\mathcal{E} } \Langle f(\{u\}), [\partial_n v] \Rangle_e \\  
&-\sum_{e\in\mathcal{E}^I } \Langle f^\prime(\{u\}) \{\partial_n u\}, [v] \Rangle_e.
\end{align*}

DG methods for the Allen-Cahn equation were studied in \cite{Feng_Li15,guo_xu18,gjx16,xia2009application}
and for the Cahn-Hilliard equation in \cite{Feng_Karakashian07,xia2007local, KSS09, AKW15, FLX16a, song2017}.
We also note that a weak Galerkin method was recently proposed and analyzed for the Cahn-Hilliard equation in 
\cite{wzzz19}. 
 
\subsubsection{Isogeometric Analysis}\label{sec-4.2.4}

In recent years, the technique of isogeometric analysis has been developed with the aim at integrating 
traditional numerical discretization like the finite element methods with computer-aided geometric design tools. 
Isogeometric analysis has been studied to many application domain including the discretization of
phase field models.  The first study in the latter direction is given in 
 \cite{gomez2008isogeometric} for  Cahn-Hilliard equations. 
 Subsequent studies include applications to phase field models for topology optimization \cite{dbh2012},
phase field models of
 brittle fracture \cite{borden2014higher}, high order equations on surfaces \cite{bdq2015}, Darcy flows 
 \cite{da2018}, chemotaxis \cite{moure2018three} 
 and electrochemical reactions \cite{zxsg2016}. 
 Convergence studies have also been carried out, see for example \cite{kmb2016}. While
 many works on isogeometric analysis are based on the Galerkin approach, there are
 also works that adopted other formulations like collocation method \cite{sbs2015} and discontinuous Galerkin method 
 \cite{zxg15}.

 Isogeometric analysis adopts $C^1$ NURBS-based finite dimensional spaces so that conforming finite element framework
 can be  utilized for spatially high order differential equations.  For phase field models with regularized
 solutions, the high order approximations can be effective using NURBS basis as in the case of spectral method 
 and high order DG methods.   For a recent review and additional references, we refer to
 \cite{gomez2018computational}.

\subsection{Spatial Mixed Discretization}\label{sec-4.3}
To overcome the difficulty caused by discretization of the biharmonic operator in the Cahn-Hilliard 
equation, especially in $3$D case,  one popular approach is to discretize the Cahn-Hilliard by 
mixed finite element or mixed DG methods, which are based on rewriting \eqref{cahn_hilliard} 
as a system of two second order PDEs given by
\begin{align}\label{CH_mixed_1}
u_t+\Delta \mu &=0,\\
\mu -\varepsilon \Delta u + \frac{1}{\varepsilon} f(u) &=0, \label{CH_mixed_2}
\end{align}  
and discretizing the above system either by finite element methods or by DG methods. 

A semi-discrete approximation can be formulated as: finding $(u^h,\mu^h)$ such that for any 
$(v^h, q^h)$, it holds
\begin{align}\label{CH_mixed_3}
(u^h_t, v^h) - a(\mu^h, v^h) & =0,\\
(\mu^h, q^h) + \varepsilon a( u^h, q^h) + \frac{1}{\varepsilon} (f(u^h),q^h) &=0.  \label{CH_mixed_4}
\end{align}  

One of the early works on the fully-discrete conforming finite element approximations to 
the mixed weak formulation of \eqref{CH_mixed_1}--\eqref{CH_mixed_2} was given 
in \cite{copetti_elliott92} which has also noted an even
earlier work by Du in an unpublished preprint.  The latter adopted the same modified Crank-Nicolson in time
discretization as in \cite{du_nicolaides91} (also see \cite{du2011finite} for similar discussions on
solving the Cahn-Hillard equation on a sphere), while the former took the backward Euler in time discretization.
Mixed finite element methods for the Cahn-Hilliard equation were also studied in 
\cite{Feng_Prohl04a, Feng_Prohl05a} with the aim of deriving refined estimates in the sharp interface limit. Mixed least square finite element was studied in \cite{DGT1996}.
Mixed DG methods were developed in \cite{KSS09, AKW15, FLX16a}.
For Cahn-Hilliard equations on evolving surfaces, mixed finite element methods were
formulated and analyzed in \cite{elliott2015evolving}.
Local DG methods were also studied for Allen-Cahn and Cahn-Hilliard 
equations \cite{GLY2016,xia2007local,xia2009application}.

%%%%%%%%
\subsection{Implementations and Advantages of High Order Methods}\label{sec-4.4}

Given the smooth phase field functions, it is expected that high order methods such as high order 
DG and spectral methods can yield more competitive numerical algorithms. For instance, errors of spectral 
approximations can be  reduced exponentially as the number of spectral basis functions grow,  
which is known in the literature as spectral accuracy. Based on such observations, Fourier spectral methods,  
for example.  have become particularly effective for phase field simulations in a periodic cell. 

Since spatial differential operators subject to periodic conditions are  diagonalizable using Fourier modes, 
the coupling of spectral spatial discretization with semi-implicit in time  schemes, as well as ETD 
and IF exponential type time integration schemes, becomes particularly popular. The treatment of nonlinearity
can also be greatly simplified with a collocation implementation or pseudo-spectral approximation
\cite{chen1998applications,cheng2016second}.

Taking, for example, the Fourier spectral method for the Allen-Cahn model over a periodic domain 
$[-\pi, \pi]^d$, we let $k$ denote a $d$-dimensional index with $|k|$ denoting the $\ell^\infty$ norm 
and consider a Fourier approximation given by 
$$u^{n}_k(x) =\sum_{|k|\le K}\hat{u}^n_k e^{ikx} \in \hat{B}_K=\text{span}\{e^{ikx}\}_{|k|\le K}.$$
We can use the first order stabilized  ETD-RK1 scheme \eqref{td_ETD1} to get \cite{du04etd,du05etd}
$$
 \hat{u}^{n+1}_k  = 
 e^{-\tau k\cdot k  \hat{u}^{n}_k  
 - \tau \alpha \varepsilon^{-2} }
 +\frac{1-e^{-\tau k\cdot k  - \tau \alpha \varepsilon^{-2} }}
 {\tau k\cdot k  \varepsilon^2 + \tau \alpha  }
   (\widehat{P_K[f(u^n_k)]}_k)  - \alpha  \hat{u}^n_k),
   $$
where $\tau$ is a step-size, $\alpha>0$ is a stabilizing constant that is often taken to be larger than
the one half of the Lipshitz constant of the nonlinear term $f$, and $P_K$ denotes the projection to $\hat{B}_K$. 
 
Replacing  $e^{\tau k\cdot k  + \tau \alpha \varepsilon^{-2} }$  by $1 + \tau k\cdot k  + \tau \alpha \varepsilon^{-2}$, we  recover
$$
(1+\tau k\cdot k  + \tau \alpha \varepsilon^{-2} ) \hat{u}^{n+1}_k  = (1+\tau \alpha \varepsilon^{-2} ) \hat{u}^n_k  - \tau\varepsilon^{-2} \widehat{P_K[f(u^n_k)]}_k\,.
$$

Computationally, Fourier spectral methods can be implemented via FFT so that the complexity per step is
$O(K\log K)$. Pseudo-spectral or spectral collocation methods are often adopted to treat the nonlinear term
instead of the spectral projection, that is, $f(u^n_k)$ is sampled at the discrete spatial grid points, then the FFT
is applied to obtained a representation in the Fourier space. 

The low computational cost and high spectral accuracy make spectral methods, when applicable, highly attractive for 
phase field simulations. An interesting question pertains to their performance 
with $\varepsilon <<1$, that is, the sharp interface limit.  In \cite{chen1998applications}, it has been shown. for
a typical benchmark setting concerning  a single spherical droplet, Fourier approximations can be much more
effective than low order finite difference even as $\varepsilon$ gets significantly reduced.  In other words, 
high order methods are competitive even in the sharp interface limit.  The latter is intriguing as the sharp interface
limit of phase field variables are generically discontinuous functions that often do not share high order 
approximations due to the lack of regularity.  A justification was given in 
\cite{dz09numerical} when mesh and time step refinement are considered as $\varepsilon$ is getting smaller.
If the sharp interface limit is viewed as the quantity of interests to be sought after, then the typical
error expected is of the form 
$$O(\varepsilon^\beta)+O(\left(\frac{\tau}{\varepsilon}\right)^\gamma)+O( 
\left(\frac{h}{\varepsilon}\right)^\delta)$$
for a finite difference or finite element methods
and
$$O(\varepsilon^\beta)+O(
\left(\frac{\tau}{\varepsilon}\right)^\gamma)+O( e^{-c K \varepsilon } )$$
 for spectral methods
where $O(\varepsilon^\beta)$ accounts for the diffuse interface approximation error while the other 
two terms are due to time  and spatial discretization.
Thus as $\varepsilon \to 0$, to have a balanced total error, asymptotically the need for refinement 
of spectral approximation is insignificant in comparison with a low order spatial approximation.

\section{Convergence Theories of Fully Discrete Numerical Methods}\label{sec-5}
%	It is clear that neither space- nor time-semidiscrete numerical schemes are 
%	implementable on computers for numerical simulations, only fully discrete schemes, 
%	which are discrete in both space and time variables, are implementable on computers. 
	In this section we shall first present a few exemplary fully discrete numerical schemes for 
	phase field models by combining the time-stepping schemes of Section \ref{sec-3} and the spatial 
	discretization methods of Section \ref{sec-4},  we then discuss the convergence results of these fully 
	discrete numerical methods. Below we divide those results into two groups: the first group addresses 
	convergence and error estimates for a fixed diffusive interfacial width $\varepsilon$, the second group concerns 
	the convergence and error estimates for the numerical solutions as well as the numerical interfaces 
	as $\varepsilon,\tau, h\to 0$ 

%\begin{enumerate}
%        \item Convergence and error estimates for a fixed interaction length $\epsilon$.
%        \item Convergence to sharp interface limits.
%\end{enumerate}

\subsection{Construction of Fully Discrete Numerical Schemes}

As for most evolution equations, a fully discrete numerical method for phase field models can be 
easily constructed by any combination of a time-stepping scheme of Section \ref{sec-4} and a 
spatial discretization method of Section \ref{sec-5}. Such a construction is also known as 
the {\em method of lines} in the literature (cf. \cite{schiesser91, reddy_trefethen92} and the references therein).
Below we shall use the Allen-Cahn equation \eqref{phase_field_5a} and the Cahn-Hilliard equation 
\eqref{cahn_hilliard} to demonstrate the procedure. 

We start with the spatial semi-discrete finite element method for the Allen-Cahn equation (with the Neumann boundary
condition),  which is given by \eqref{AC_GalerkinApprox} with \eqref{AC_FEM}, since $V_h$ is finite-dimensional, it is 
easy to see that \eqref{AC_GalerkinApprox} with \eqref{AC_FEM} is a system of ODEs, which must be complemented 
by the initial condition 
\[
\bigl( u_h(0), v_h \bigr) = \bigl( u_0, v_h \bigr)\qquad \forall v_h\in V_h. 
\]

Applying the convex splitting scheme to \eqref{AC_GalerkinApprox} then leads to the following 
fully discrete convex splitting finite element method: find $u^{n+1}_h\in V_h$ such that 

\begin{align}\label{AC_GalerkinApprox_CS}
\bigl(d_t u_h^{n+1}, v_h\bigr) + a_h\big( u_h^{n+1}, v_h\bigr) + \frac{1}{\varepsilon^2} \Bigl((u_h^{n+1})^3, v_h\Big)
= \frac{1}{\varepsilon^2} \bigl(u_h^n, v_h\big) %\qquad \forall v_h\in V_h,
\end{align}  
for all $v_h\in V_h$. Where $d_tu_h^{n+1}:= (u_h^{n+1} -u_h^n)/\tau$. As expected, it is easy to check that
the above method is unconditionally (in $h$ and $\tau$) and uniformly (in $\varepsilon$) energy stable (cf.   \cite{shen2010,Feng_Li15}).

Similarly, applying the implicit-explicit (IMEX) scheme to \eqref{AC_GalerkinApprox} gives the following 
fully discrete IMEX finite element method: find $u^{n+1}_h\in V_h$ such that

\begin{align}\label{AC_GalerkinApprox_IMEX}
\bigl(d_t u_h^{n+1}, v_h\bigr) + a_h\big( u_h^{n+1}, v_h\bigr)   
= \frac{1}{\varepsilon^2} \Bigl(u_h^n -(u_h^n)^3, v_h\Big) ,\quad \forall  v_h\in V_h. 
\end{align}  
 As expected, this IMEX finite element method is not unconditional stable. To strengthen its stability, the following  popular stabilized method was proposed in the literature 
\cite{xu_tang06, LQT16, FTY15, li2016, shen2010}: find $u^{n+1}_h\in V_h$ such that
\begin{align}\label{AC_GalerkinApprox_IMEX_stab}
\bigl(d_t u_h^{n+1}, v_h\bigr) + a_h\big( u_h^{n+1}, v_h\bigr) +S \bigl( u_h^{n+1}-u_h^n, v_h\bigr)   
= \frac{1}{\varepsilon^2} \Bigl(u_h^n -(u_h^n)^3, v_h\Big)
\end{align}  
for all $v_h\in V_h$. Where $S$ is an under-determined constant or linear operator. 

Now adopting the finite element method for the spatial discretization in the IEQ scheme \eqref{CH_IEQ} 
then leads to the following fully discrete IEQ finite element method for the Cahn-Hilliard equation \eqref{cahn_hilliard}: find $\bigl( u^{n+1}_h, \mu^{n+1}_h, q^{n+1}_h\bigr)\in [V_h]^3$ such that
 %\begin{subequations}\label{CH_IEQ_FE}
 \begin{align*}
 \bigl( u^{n+1}_h - u^n_h, v_h \bigr) &= -\tau a_h\bigl( \mu^{n+1}_h, v_h \bigr) \quad \forall v_h\in V_h,\\
 \bigl(\mu^{n+1}_h, \lambda_h \bigr) &= \varepsilon a_h\bigl( u^{n+1}_h, \lambda_h \bigr) 
 + \frac{1}{\varepsilon} \Bigl(  G(u^n_h)  q^{n+1}_h,\lambda_h \Bigr) \quad \forall \lambda_h\in V_h,\\
 \bigl(q^{n+1}_h -q^n_h, p_h\bigr) & =  \frac12 \Bigl( G(u^n_h) \bigr(u^{n+1}_h-u^n_h\bigr),  p_h \Bigr)  
  \quad \forall p_h\in V_h, 
 \end{align*}
 %\end{subequations}
where $a_h(\cdot,\cdot)$ denotes the discrete bilinear form for $-\Delta$ as defined in \eqref{AC_FEM}.

Finally, using the finite element method for the spatial discretization in the SAV scheme \eqref{CH_SAV} 
one immediately obtains the following fully discrete SAV finite element method for the Cahn-Hilliard equation
 \eqref{cahn_hilliard}: find $\bigl( u^{n+1}_h, \mu^{n+1}_h, r^{n+1}_h\bigr)$ $\in [V_h]^2\times \mathbf{R}$ such that
%$$%\begin{subequations}\label{CH_SAV_FE}
\begin{align*}
\bigl( u^{n+1}_h - u^n_h, v_h \bigr) &= -\tau a_h\bigl( \mu^{n+1}_h, v_h \bigr) \quad \forall v_h\in V_h,\\
\bigl(\mu^{n+1}_h, \lambda_h \bigr) &= \varepsilon a_h\bigl( u^{n+1}_h, \lambda_h \bigr) 
+ \frac{1}{\varepsilon} \Bigl(  H(u^n_h), \lambda_h \Bigr)   r^{n+1}_h \quad \forall \lambda_h\in V_h,\\
r^{n+1}_h -r^n_h & =  \frac{1}{2} \Bigl( H(u^n_h), u^{n+1}_h-u^n_h \Bigr), 
\end{align*}
%$$%\end{subequations}
where $a_h(\cdot,\cdot)$ is the same as above.  

\subsection{Types of Convergence and a Priori Error Estimates}

As for many numerical methods, three important theoretical issues to be addressed for above fully 
discrete numerical methods are {\em stability, convergence} 
and {\em rates of convergence}. As the stability issue has already been considered when constructing 
the time-stepping schemes, and it is not difficult to show that the similar stability properties 
also hold for the corresponding fully discrete schemes, especially when Galerkin spatial discretization 
methods are used, in this section we only focus on the issues related to  
convergence and rates of convergence (i.e., error estimates) of fully discrete numerical methods for 
phase field models in two cases: (i) $\varepsilon>0$ is fixed and $h, \tau\to 0$; 
(ii) $\varepsilon, h, \tau\to 0$. It should be emphasized that the stability is only a necessary 
criterion for screening ``good" numerical methods, however, it does not guarantee 
the convergence, as a result, it is important to examine the convergence (and rates 
of convergence) for stable methods, see \cite{XLWB18} for further discussions in this direction. 

Let $u^\varepsilon$ denote the PDE solution of a underlying phase field model and $u^\varepsilon_{h,\tau}$
denote a fully discrete numerical solution.  The goal of error estimates for a fixed $\varepsilon>0$ 
is to derive the following type of error bounds:
\begin{equation}\label{error_estimate}
\|u^\varepsilon- u^\varepsilon_{h,\tau}\|_* \leq C_1(\varepsilon, u^\varepsilon) \tau^{\ell_1} 
+C_2(\varepsilon,u^\varepsilon) h^{\ell_2}
\end{equation}
for some $\varepsilon$- and $u^\varepsilon$-dependent positive constants $C_1$ and $C_2$ 
and positive numbers (mostly integers) $\ell_1$ and $\ell_2$. Where $\|\cdot\|_*$ denotes 
a (function space) norm which is problem-dependent. 

The precise dependence of $C_1$ and $C_2$ on $\varepsilon$ are usually complicated and difficult to 
obtain, but they are expected to grow in $\frac{1}{\varepsilon}$ as $\varepsilon\to 0$, this is because 
$C_1$ and $C_2$ depend on high order ($\geq 2$) space and time derivatives of $u^\varepsilon$ and those 
derivatives can be proved to grow {\em polynomially} in $\frac{1}{\varepsilon}$ as $\varepsilon\to 0$
(cf. \cite{Feng_Prohl03b,Feng_Prohl05a,Feng_Prohl04c,FHL06,Feng06,abels_lengeler14, abels15, ALS17, abels_liu18}). 

Here we distinguish two cases based on how $C_1$ and $C_2$ depend on $\frac{1}{\varepsilon}$.
If at least one of $C_1$ and $C_2$ grow {\em exponentially} in $\frac{1}{\varepsilon}$,
we then call \eqref{error_estimate} a {\em coarse error estimate}. On the other hand, if both $C_1$ and $C_2$ 
grow {\em polynomially} in $\frac{1}{\varepsilon}$, the estimate is called a {\em fine error estimate}. 
Obviously, in both cases, for a fixed $\varepsilon>0$, \eqref{error_estimate} implies the convergence 
\[
\lim_{h,\tau\to 0} \|u^\varepsilon- u^\varepsilon_{h,\tau}\|_*=0,
\]
that is, $u^\varepsilon_{h,\tau}$ converges to $u^\varepsilon$ in the $\|\cdot\|_*$-norm. Hence, 
when $\varepsilon>0$ is fixed, both types of coarse and fine error estimates imply the convergence 
of the numerical solution as $h,\tau\to 0$. 

On the other hand,  the situation is very different if the convergence of the numerical interface 
to the sharp interface of the underlying geometric moving interface problem is a main concern. 
In that case, one must study the limiting behaviors of the numerical solution $u^\varepsilon_{h,\tau}$
and the numerical interface 
\[
\Gamma^{\varepsilon,h,\tau}_t:=\bigl\{ x\in \Omega;\, u^\varepsilon_{h,\tau}(x,t)=0 \bigr\}
\]
as $\varepsilon\to 0$.  
First, it is easy to see that the coarse error estimates will become useless because they
fail to provide any useful information about the convergence. To see the point, we notice that
$\tau^{\ell_1}$ and $h^{\ell_1}$ only decrease in polynomial orders but $C_1$ and/or $C_2$ grow in 
exponential orders, hence, the limit of the right-hand side of \eqref{error_estimate} is $+\infty$
as $\varepsilon, h, \tau\to 0$. 
Second, on contrary the fine error estimates are still valuable, this is because $C_1$ and $C_2$ grow in 
polynomially orders, $h$ and $\tau$ can be chosen as powers of $\varepsilon$ so that 
the right-hand side of \eqref{error_estimate} is guaranteed to converge to $0$ as $\varepsilon\to 0$. 
Let $u^0$ denote the solution of the underlying moving sharp interface problem (which is defined by its 
level set formulation), by the triangle inequality we have 
\[
\|u^0- u^\varepsilon_{h,\tau}\|_\dagger \leq \|u^0- u^\varepsilon \|_\dagger
+ \|u^\varepsilon- u^\varepsilon_{h,\tau}\|_\dagger,
\]
where $\|\cdot\|_\dagger$ stands for another (function space) norm which is a topology used to measure 
the convergence of $u^\varepsilon$ to $u^0$ as described in Section \ref{sec-2}. We note that 
$\|\cdot\|_\dagger$ is usually a weaker norm compared to the $\|\cdot\|_*$-norm. Thus, combining the PDE 
convergence result and the above numerical convergence result we get 
\[
\lim_{\varepsilon\to 0} \|u^0-u^\varepsilon_{h,\tau}\|_\dagger =0,
\]
provided that $h$ and $\tau$ are chosen as appropriate powers of $\varepsilon$. Moreover, it turns out that  
\cite{Feng_Prohl03b, Feng_Prohl04c, Feng_Prohl05a,Feng_Li15,FLX16a} such a 
convergence result is good enough to infer the following convergence of the numerical interface:
\[
\lim_{\varepsilon\to 0} \mbox{dist}_H \bigl(\Gamma_t, \Gamma^{\varepsilon,h,\tau}_t \bigr) =0,
\]
as well as its rates of convergence in powers of $\varepsilon$ if $\Gamma_t$ is sufficiently smooth.
Where $\mbox{dist}_H$ denotes the Hausdorff distance between two sets. 

We shall focus on coarse error estimates in the next subsection and discuss fine error estimates 
and the convergence of numerical interfaces in the subsequent subsection.

\subsection{Coarse Error Estimates for a Fixed Value $\varepsilon>0$}

For a fixed $\varepsilon>0$, all phase field equations are either semilinear or quasilinear  
parabolic PDEs, the desired error estimates for their fully discrete (and spatial semi-discrete) methods can be 
obtained by following the standard perturbation procedure which consists of three main steps: (i) decomposing 
the global error with help of an elliptic projection of the exact PDE solution; (ii) using the energy
method to obtain an error equation/inequality for the error between the numerical solution and the elliptic 
projection; (iii) completing the error estimate using a discrete Gronwall's inequality and the triangle 
inequality.  Below we shall demonstrate this perturbation procedure using scheme \eqref{AC_GalerkinApprox_CS}. 

Let 
\[
R^{n+1}:= u_t (t_{n+1}) - d_t u^{n+1} 
\]
denote the truncation error of the backward difference operator $d_t$ which is well known to be 
of $O(\tau)$ order provided that $u_{tt} \in L^2((0,T); H^{-1}(\Omega))$. Subtracting 
\eqref{AC_GalerkinApprox_CS} from the weak formulation of \eqref{phase_field_5a} immediately yields 
the following error equation:
\begin{align}\label{error_eq_AC0}
	\bigl(d_t e_h^{n+1}, v_h\bigr) + a_h\bigl(e_h^{n+1}, v_h\bigr) 
	&+ \frac{1}{\varepsilon^2} \bigl( (u(t_{n+1}))^3-(u_h^{n+1})^3, v_h \bigr) \\
	&= \frac{1}{\varepsilon^2} \bigl( e_h^{n+1} + \tau d_t u_h^{n+1}, v_h\bigr)  + R^{n+1} \quad \forall v_h\in V_h, \nonumber
\end{align}
where $e_h^n:= u(t_n) - u_h^n$ denotes the global error at $t_n$. Note that the super-index $\varepsilon$ is 
suppressed on all functions. 

Now, let $P_h: H^1(\Omega)\to V_h$ denote the standard elliptic projection operator associated with
the discrete bilinear form $a_h(\cdot,\cdot)$ and set  
\[
e_h^n:= \eta_h^n +\xi_h^n; \quad \eta_h^n:= u(t_n)-P_h u(t_n), \, \,
\xi_h^n:= P_h u(t_n) -u_h^n. 
\]
Then \eqref{error_eq_AC0} implies that 
\begin{align} \label{error_eq_AC1}
\bigl(d_t \xi_h^{n+1}, v_h\bigr) + a_h\bigl(\xi_h^{n+1}, v_h\bigr) 
&+\frac{1}{\varepsilon^2} \bigl( (P_hu(t_{n+1}))^3-(u_h^{n+1})^3 -\xi_h^{n+1} , v_h \bigr) \\
&= \frac{1}{\varepsilon^2} \bigl( \eta_h^{n+1}+ \tau d_t u_h^{n+1}, v_h\bigr) 
+ R^{n+1} -\bigl( d_t \eta_h^{n+1}, v_h\bigr) \nonumber \\
&\quad -\frac{1}{\varepsilon^2} \bigl( (u(t_{n+1}))^3-(P_h u_h^{n+1})^3, v_h \bigr) 
\quad \forall v_h\in V_h. \nonumber
\end{align}
Setting $v_h=\xi_h^{n+1}$ and using the monotonicity of cubic power function and the 
stability estimates of the numerical solution $u_h^\varepsilon$ (not shown here) as well as
the approximation properties of $P_h$ we get 
\begin{align}\label{error_eq_AC}
\frac12 d_t \|\xi_h^{n+1}\|^2_{L^2} + a_h\bigl(\xi_h^{n+1}, \xi_h^{n+1}\bigr) 
&\leq \frac{1}{\varepsilon^2} \| \xi_h^{n+1}\|^2_{L^2}  + c_1 \tau^2   + c_2 h^{r+1}
\end{align}
for some positive constants $c_i=c_i(\varepsilon) (i=1,2)$.  Where $r$, a positive integer, 
denotes the order of the underlying finite element method 
(i.e., the degree of polynomial shape functions). Here we have also used the fact that
\[
\|\eta_h^{n+1}\|_{L^2} \leq C h^{r+1}.
\]

Finally, applying a discrete Gronwall's inequality to \eqref{error_eq_AC} yields
\begin{align}\label{error_eq_AC2}
\|\xi_h^m\|^2_{L^2} + \tau \sum_{n=0}^{m-1} a_h\bigl(\xi_h^{n+1}, \xi_h^{n+1}\bigr) 
\leq C\exp \bigl(\varepsilon^{-2} \bigr) \bigl( c_1 \tau^2   + c_2 h^{r+1} \bigr)
\end{align}
for some positive constant $C=C(T)$, which and the triangle inequality infer 
\begin{align}\label{L2_error}
\|e_h^m\|^2_{L^2} &\leq C\exp \bigl(\varepsilon^{-2} \bigr) \bigl( c_1 \tau^2   + c_2 h^{r+1} \bigr),\\
\tau \sum_{n=0}^{m-1} a_h\bigl(e_h^{n+1}, e_h^{n+1}\bigr) &\leq C\exp \bigl(\varepsilon^{-2} \bigr) \bigl( c_1 \tau^2  
 + c_2 h^r \bigr) \label{H1_error}
\end{align}
for all $1\leq m\leq M$. 

We remark that error estimates \eqref{L2_error} and \eqref{H1_error} are two coarse error estimates for the 
fully discrete implicit Euler finite element method for the Allen-Cahn equation. It turns out that
they are typical error estimates for all fully discrete numerical methods for various phase field models, which 
can be abstractly written into the form \eqref{error_estimate}. 
Although for different models, the norms used to measure errors and 
the orders of errors (indicated by $\ell_1$ and $\ell_2$ in \eqref{error_estimate} ) may be different, 
the nature of the exponential dependence on $\frac{1}{\varepsilon}$ of the constants $C_1$ and $C_2$
remains the same for all these methods and models when the above standard perturbation procedure 
is used to derive the error estimates, see  
 \cite{elliott_french87, elliott_french89, du_nicolaides91, barrett_blowey99, barrett_blowey_garcke99, xu_tang06,
 	FTY15, liu_shen03,du2004phase, FHL06, Feng06, Feng_Karakashian07, wise2009, LLRV09, KSS09, shen2010, AKW15, 
 	li2016, diegel2017convergence, yang16, song2017, shen2018}.

%%%%%%%
\subsection{Fine Error Estimates and Convergence of Numerical Interfaces as $\varepsilon, h,\tau\to 0$}

As explained earlier, the coarse error estimates obtained in the previous subsection are not 
useful to study the convergence as $\varepsilon\to 0$.  A couple natural questions arise: 
(i) Are the error estimates sharp in terms of $\varepsilon$?  (ii) If not, what are sharp error 
estimates and how to establish such estimates? As all numerical experiments (cf. \cite{dz09numerical, PDG99}
and the references therein) indicate that the 
above coarse error estimates are not sharp in terms of $\varepsilon$ (otherwise, the 
phase field methodology would be an impractical method to compute moving interface problems!), 
then the focus is on addressing the second question.  The good news is that there are 
a couple successful numerical analysis techniques which have been developed in the last thirty years, 
the bad news is that they are the only known techniques for deriving so-called fine error estimates, and 
the applicability of both techniques is restricted although one technique fares better than the other. 
Below we shall again use scheme \eqref{AC_GalerkinApprox_CS} for the Allen-Cahn equation 
as an example to explain the ideas of both techniques.

The first technique, which was developed in \cite{Nochetto_Verdi97,nochetto_paolini_verdi94},  
is to use the discrete maximum principle to derive the $L^\infty$-norm {\em fine} error estimate 
for $u^\varepsilon-u^\varepsilon_{h,\tau}$ for the implicit Euler $P_1$-conforming finite element  
(i.e., $r=1$), such an estimate in turn leads to the convergence proof of the numerical 
interface $\Gamma_t^{\varepsilon,h,\tau}$ to its underlying sharp interface limit $\Gamma_t$.
However, this maximum principle technique only applies to the Allen-Cahn equation which 
does possesses a maximum principle and only to $P_1$-conforming finite element method in order to 
to ensure a discrete maximum principle. The technique does not work for high order phase 
field model such as the Cahn-Hilliard equation nor for phase field systems which do not have a 
maximum principle.  Moreover, it does not apply to higher order finite element methods either 
because higher order finite element methods do not have a discrete maximum principle in general even 
the underlying PDE (such as the Allen-Cahn equation) does. 

Before introducing the second technique, we first have a closer look at the derivation of
the coarse error estimate given in the previous subsection to find out the guilty part/step
which leads to the exponential dependence of the error constants on $\frac{1}{\varepsilon}$. 
First, since the error equation \eqref{error_eq_AC1} is an identity, no error is introduced there.
Second, after the inequality \eqref{error_eq_AC} is reached, noticing that the coefficient 
$\frac{1}{\varepsilon^2}$ in the first term on the right-hand side,  then it is too late to 
improve the estimate because \eqref{error_eq_AC} and the Gronwall's inequality 
certainly lead to a coarse error estimate!  This simple observation suggests that 
the guilty step is from the error equation \eqref{error_eq_AC1} to the error inequality 
\eqref{error_eq_AC}, in other words, the transition from equality to inequality is too loose
which makes the right-hand side of \eqref{error_eq_AC} becomes a too large bound for its 
left-hand side. To improve the coarse error estimate, we must refine this step and to 
obtain a tight upper bound for the left-hand side of \eqref{error_eq_AC}. It turns out
that this is exactly what needs to be done to derive a fine error estimate for the Allen-Cahn equation,
and a similar situation occurs for several other phase field model. 

Moreover, since the right-hand side of \eqref{error_eq_AC1} only involve the truncation error 
and the projection error, they contribute $c_1\tau^2$ and $c_2 h^{r+1}$ terms in \eqref{error_eq_AC},
which is good.  Hence, the desired improvement (if possible) must comes from bounding 
the last term on the left-hand side of \eqref{error_eq_AC1}. To obtain \eqref{error_eq_AC}
we bounded these term as two separate terms, that inevitably produces the first term 
on the right-hand side of \eqref{error_eq_AC}. A key idea of the second technique of the fine 
error analysis, which was developed in \cite{Feng_Prohl03b,Feng_Prohl04c,Feng_Prohl05a,Feng_Li15,
FLX16a}, is not to separate that term, instead, to estimate it together with the proceeding term.  
We demonstrate the steps of this techniques below. 

Let $f(u):=u^3-u$, then $f'(u)=3u^2-1$. The last term on the left-hand side of \eqref{error_eq_AC1} with 
$v_h=\xi_h^{n+1}$ can be 
written as 
\begin{align*}  
% \frac{1}{\varepsilon^2} \bigl( (P_hu(t_{n+1}))^3 &-(u_h^{n+1})^3 -\xi_h^{n+1} , \xi_h^{n+1} \bigr)  \\
% & =\frac{1}{\varepsilon^2} \bigl( f(P_h u(t_{n+1})) -f(u_h^{n+1}) , \xi_h^{n+1} \bigr) \\
 \frac{1}{\varepsilon^2} \Bigl( f'(P_h u(t_{n+1})) \xi_h^{n+1}, &\,\xi_h^{n+1} \Bigr) 
     + \frac{1}{\varepsilon^2} \| \xi_h^{n+1} \|_{L^4}^4
    - \frac{3}{\varepsilon^2} \Bigl(P_h u(t_{n+1}), (\xi_h^{n+1})^3   \Bigr) \\
 &\geq  \frac{1}{\varepsilon^2} \Bigl( f'(P_h u(t_{n+1})) \xi_h^{n+1}, \xi_h^{n+1} \Bigr) 
 + \frac{1}{\varepsilon^2} \| \xi_h^{n+1} \|_{L^4}^4 \\
 &\qquad\quad  - \frac{3}{\varepsilon^2} \| P_h u(t_{n+1})\|_{L^\infty} \|\xi_h^{n+1}\|_{L^3}^3.
\end{align*}
Substituting the above inequality into \eqref{error_eq_AC1} yields
\begin{align}\label{fine_error_eqn}
&\frac12 d_t \|\xi_h^{n+1}\|^2_{L^2} + \varepsilon^2\Bigl[ a_h\bigl(\xi_h^{n+1}, \xi_h^{n+1}\bigr) 
+ \frac{1}{\varepsilon^2} \Bigl( f'(P_h u(t_{n+1})) \xi_h^{n+1}, \xi_h^{n+1} \Bigr)\Bigr]  \\
&\,\,  + (1-\varepsilon^2)  \Bigl[ a_h\bigl(\xi_h^{n+1}, \xi_h^{n+1}\bigr) 
+ \frac{1}{\varepsilon^2} \Bigl( f'(P_h u(t_{n+1})) \xi_h^{n+1}, \xi_h^{n+1} \Bigr)\Bigr]  
+ \frac{1}{\varepsilon^2} \| \xi_h^{n+1} \|_{L^4}^4   \nonumber \\
& \hskip 1.52in   \leq \frac{C}{\varepsilon^2}  \|\xi_h^{n+1}\|_{L^3}^3   + c_1 \tau^2   + c_2 h^{r+1}, \nonumber
\end{align}
which replaces \eqref{error_eq_AC}. One key step of the technique is to show  
that there exists an $\varepsilon$-independent positive constant $\widehat{C}_0$ such that  
\begin{align}\label{fine_error_eqn2}
a_h\bigl(\xi_h^{n+1}, \xi_h^{n+1}\bigr) 
+ \frac{1}{\varepsilon^2} \Bigl( f'(P_h u(t_{n+1})) \xi_h^{n+1}, \xi_h^{n+1} \Bigr) 
\geq -\widehat{C}_0\,\|\xi_h^{n+1}\|^2_{L^2},
\end{align}
which is a corollary of the following PDE spectral estimate result \cite{deMottoni_schatzman95,chen94}:
$$%\begin{equation} \label{discrte_spectral_estimate_AC}
\lambda_{\mbox{\small min}}:= \inf_{\varphi \in H^1} \frac{ \bigl(\nabla \varphi, \nabla \varphi \bigr)  
	+ \frac{1}{\varepsilon^2} \Bigl( f'(u^\varepsilon) \varphi, \varphi \Bigr) }{ \|\varphi\|_{L^2} } \geq -C_0
$$%\end{equation}
for some positive constant $C_0$ and $0<\varepsilon <<1$. Notice that $\lambda_{\mbox{\small min}}$ 
is the principal eigenvalue of the linearized Allen-Cahn operator
%\begin{equation}\label{PDE_spectral_estimate_AC}
$$L_{\mbox{\small AC}} (\varphi):= -\Delta \varphi + \frac{1}{\varepsilon^2} f'(u^\varepsilon) \varphi.  $$
%\end{equation} 

Another key step of the technique is to bound the first term on the right-hand side of 
\eqref{fine_error_eqn} as follows: 
\begin{align*}%\label{fine_error_eqn3} 
\frac{C}{\varepsilon^2} \|\xi_h^{n+1}\|_{L^3}^3 
%&\leq \frac{C}{\varepsilon^2} \Bigl( 
%\|\nabla \xi_h^{n+1}\|_{L^2}^{\frac{d}{2}}   \| \xi_h^{n+1}\|_{L^2}^{\frac{6-d}{2} }
%+  \|\xi_h^{n+1}\|_{L^2}^3 \Bigr) \\
%&\leq  \frac{\varepsilon^2}{2} \|\nabla \xi_h^{n+1}\|_{L^2}^2 
%+ C \varepsilon^{\frac{-2(4+d)}{4-d}} \|\xi_h^{n+1}\|_{L^2}^{\frac{2(6-d)}{4-d}} 
%+ \frac{C}{\varepsilon^2} \|\xi_h^{n+1}\|_{L^2}^{1+2} \nonumber \\
&\leq  \frac{\varepsilon^2}{2} \|\nabla \xi_h^{n+1}\|_{L^2}^2 
+ \widehat{C}_0 \varepsilon^2 \|\xi_h^{n+1}\|_{L^2}^2 
+  C \varepsilon^{\frac{-2(4+d)}{4-d}} \|\xi_h^{n+1}\|_{L^2}^{\frac{2(6-d)}{4-d}}. %\nonumber
\end{align*}

Combining  the above with \eqref{fine_error_eqn} and  \eqref{fine_error_eqn2}  yields
\begin{align}\label{fine_error_eqn4}
&d_t \|\xi_h^{n+1}\|^2_{L^2} + \varepsilon^2a_h\bigl(\xi_h^{n+1}, \xi_h^{n+1}\bigr) 
+   \frac{2}{\varepsilon^2} \| \xi_h^{n+1} \|_{L^4}^4     \\
&\qquad  \leq C \varepsilon^{\frac{-2(4+d)}{4-d}} \|\xi_h^{n+1}\|_{L^2}^{\frac{2(6-d)}{4-d}}
+ 2\bigl(1+\widehat{C}_0 \bigr) \|\xi_h^{n+1}\|^2_{L^2}
+ c_1 \tau^2 + c_2 h^{r+1}. \nonumber
\end{align}
Notice that the exponent $\frac{2(6-d)}{4-d}> 2$ for all $d\geq 1$, the final key step of the technique is 
to obtain the desired fine error bound for $\xi_h^{n}$ (and hence for $e_h^n$) from \eqref{fine_error_eqn4}
by using a discrete generalized Grownwall's inequality. See \cite{Feng_Prohl03b} for finite element 
approximations and \cite{Feng_Li15} for discontinuous Galerkin approximations. 

From the above derivation one can see that the most important step of the second technique for fine error 
estimates is to establish a discrete spectral estimate based on the PDE spectral estimate for the linearized 
operator of the underlying PDE operator. It was proved that such a spectral estimate also holds for the 
Cahn-Hilliard equation/operator and 
for the general phase field model for the generalized Stefan problem 
\cite{bates_fife90, alikakos_fusco93, chen94}, and similar 
fine error estimates to above estimate were also established in \cite{Feng_Prohl04a,Feng_Prohl05a} 
(mixed finite element approximations), \cite{FLX16a} (mixed discontinuous Galerkin approximations),
and \cite{wu_li18, Li18} (Morley nonconforming finite element approximation) for the Cahn-Hilliard equation 
and in \cite{Feng_Prohl04c} for the general phase field model based on those spectral estimates. 
It should be emphasized that this technique is applicable to any phase field model as long as its linearized operator  
satisfies a desired spectral estimate. Hence, this technique essentially reduces a numerical problem of 
deriving fine error estimates for a phase field model into a problem of proving a PDE spectral estimate for the 
linearized phase field operator. 
Finally, we also note that a duality argument was recently proposed in \cite{Chrysafinos17} 
to derive unconditional stability for a space-time finite element method, it may have a potential to 
offer an alternative technique for fine error estimates.

\section{A Posteriori Error Estimates and Adaptive Methods}\label{sec-6}
%Topics to be covered include (not limited)
%\begin{enumerate}
%        \item Ad hoc methods ($h$-, $hp$- and $r$-adaptive methods).
%        \item Methods with convergence analysis.
%\end{enumerate}

\subsection{Spatial and Temporal Adaptivity}

Since phase field models are singularly perturbed equations which involve a small 
scale parameter $\varepsilon$ and the solutions of phase field models have distinctive profiles
in the sense that they take a close-to-constant value in each bulk region, which indicates a phase of the material or 
fluid mixture, and vary smoothly but sharply in a thin layer (called diffuse interface) of width 
$O(\varepsilon)$ between bulk regions, to resolve this kind of (solution) functions, it is necessary to 
use spatial mesh size smaller than $\varepsilon$ (sometimes much smaller) in the thin layer and it is 
preferable to use much coarser meshes in the bulk regions. This seems to provide an ideal 
situation for using adaptive grid methods for efficient numerical simulations. Moreover, 
should a uniform mesh be used for simulations, the mesh must be very fine over the whole domain, 
this then results in huge nonlinear algebraic systems to solve, which may not be feasible 
for large scale and long time simulations, especially in high dimensions. 
Furthermore, the time scale in most phase field models represents so-called fast time scale, to capture 
the dynamics of underlying physical or biological phase transition phenomena, it is necessary to use 
very fine time step size. If the physical or biological processes are slowly varying, it often takes 
a long time for them to reach equilibrium states, when that happens, long time numerical simulations 
are required and adaptive time-stepping schemes must be used to make such numerical simulations feasible \cite{li2017computationally,chen2016efficient,wodo2011computationally,zhang2012adaptive}. 

Since the method of lines is the preferable approach in practice, as explained in 
Sections \ref{sec-3}--\ref{sec-5}, after a spatial discretization is done, each phase field model then 
reduces into a system of nonlinear ordinary differential equations (ODEs), then the existing 
well-developed adaptive time-stepping schemes for ODEs can be utilized for temporal adaptivity. 
For this reason we shall not further discuss temporal adaptivity in this section and refer the 
reader to \cite{Iserles09, Lambert91}, instead, below we shall focus on discussing the ideas of various
approaches of spatial adaptivity for phase field models. 

To design spatial adaptive methods/algorithms for PDE problem, three major approaches have been  
extensively developed in the past thirty years, they are $h$-, $hp$- and $r$-adaptivity 
(see \cite{EEJ95, ainsworth_oden00, huang_russell11} for detailed expositions). The $h$-adaptivity 
starts with a relatively coarse initial mesh and then refines or coarsens the mesh successively based on a chosen 
error estimator, the sought-after adaptive method, if done correctly, should able to automatically refine 
or coarsen the mesh where a refinement or coarsening should be done, so the method successively searches 
an ``optimal" mesh and computes a satisfactory numerical solution on the mesh. The $h$-adaptivity often 
goes together with finite element and finite difference methods, and the same finite element 
or finite difference method is often employed on each mesh successively generated by the adaptive algorithm. 
That is, the degree of the finite element space does not change in the algorithm. A typical $h$-adaptive
algorithm consists of successive loops of the following sequence: 
\[
\boxed{\bf SOLVE}\to \boxed{\bf ESTIMATE}\to \boxed{\bf MARK}\to \boxed{\bf REFINE/COARSEN}.
\]
The $hp$-adaptivity is often associated with the $hp$-finite element method \cite{Schwab98}, it only differs from the $h$-adaptivity in one aspect, that is, the refinement or coarsening is not only done 
for the mesh but also for the degree of finite element spaces. Hence, the finite element spaces 
may be different on the successive meshes generated by the adaptive method. 
The $r$-adaptivity, which is significantly different from the $h$- and $hp$-adaptivity, uses a 
fixed number of mesh points in the adaptive algorithm, however, the distribution of these mesh 
points is allowed to move, hence, a new mesh is generated by a new distribution of the mesh points,
in particular, a locally fine mesh can be obtained by clustering a large (enough) number of the mesh 
points to where a fine mesh is desired. Clearly, a central issue for $r$-adaptive methods is how to 
develop an automatic strategy to move the mesh points.  
In \cite{feng2006spectral,feng2009fourier,yu2008applications,shen2009efficient}, $r$-adaptivity was developed 
for spectral approximations of several phase field models. For more strategies on $r$-adaptivity, we refer to
 \cite{huang_russell11,tan2006simple,zhang2007adaptive,shen2011spectral, hu_li_tang09,DLT08}. 

In the rest of this section, we shall only discuss $h$-adaptive finite element methods for phase field models. Indeed, there has been a lot of interest in
spatial adaptive simulations for various phase field models 
 \cite{ceniceros2007nonstiff,ceniceros2010three,du2008adaptive,JJ2018,rosam2008adaptive,provatas2005multiscale,braun1997adaptive,stogner2008approximation}.

We refer the reader to \cite{ainsworth_oden00, shen2011spectral,huang_russell11}  for detailed expositions for the other two 
spatial adaptive methodologies.

\subsection{Coarse and Fine a Posteriori Error Estimates for Phase Field Models}

As mentioned above the key step for developing an $h$-adaptive method is to design a 
``good" error indicator which can successively predict where to refine and where to coarsen the mesh. 
As expected, such an indicator should be computable and must relate to the error of numerical solutions,  
so the mesh is refined where the error is 
large and coarsened where the error is too small. Error indicators of this kind is known as 
{\em a posteriori} error estimates in the literature.

The general form of {\em a posteriori} error estimates is given by 
\begin{equation}\label{posteriori_error_estimate}
\|u^\varepsilon- u^\varepsilon_{h,\tau}\|_{\#} \leq \mathcal{E}(\varepsilon, h,\tau, u^\varepsilon_h) 
\end{equation}
for some positive functional $\mathcal{E}$. Where $\|\cdot\|_{\#}$ denotes 
a (function space) norm which is problem-dependent. Compare to the error estimate \eqref{error_estimate}, 
there is one important difference, that is, $\mathcal{E}(\varepsilon, h,\tau, u^\varepsilon_h)$ is computable 
because it depends on the numerical
solution $u^\varepsilon_h$, not the PDE solution $u^\varepsilon$; on the other hand, the error bound
in \eqref{error_estimate} is not computable because it depends on the unknown PDE solution $u^\varepsilon$.
For this very reason, \eqref{error_estimate} and \eqref{posteriori_error_estimate} are called respectively
{\em a priori} and {\em a posteriori} error estimates. 

The primary objective of {\em a posteriori} error estimates is to derive a tight error 
bound $\mathcal{E}(\varepsilon, h,\tau, u^\varepsilon_h)$ which is relatively easy to compute.
Similar to the classification of  {\em a priori} error estimates,  we also divide {\em a posteriori} 
error estimates into two groups for phase field models. The first group of estimates, called coarse estimates,
depend on $\varepsilon^{-1}$ exponentially and the second group, called fine estimates,
depend on $\varepsilon^{-1}$ polynomially. Both groups of {\em a posteriori} error estimates
have been developed for phase field models although the majority of them belong to the first group.  
Moreover, two main techniques have been used to derive the {\em a posteriori} error bound/functional  
$\mathcal{E}(\varepsilon, h,\tau, u^\varepsilon_h)$. The first technique gives rise so-called 
residual-based error bounds 
\cite{EEJ95,morin_nochetto_siebert02, Binev_Dahmen_DeVore04, nochetto_siebert_veeser09}, 
and the second one 
is a duality-based technique, called the {\em Dual Weighted Residual} method, which leads to goal-oriented 
error estimates \cite{bangerth_rannacher03}. Below we shall only focus on discussing the residual-based 
{\em a posteriori} error estimation because of its popularity and superiority for phase field models. 

The residual-based {\em a posteriori} error estimates give an error bound/functional 
$\mathcal{E}(\varepsilon, h,\tau, u^\varepsilon_h)$ which typically depends on local residuals 
of $u^\varepsilon_h$ and the jumps of the flux $\nabla u^\varepsilon_h\cdot n$ along each 
interior element edge/face. It turns out that residual-based error estimators are quite effective 
for phase field models, this is because the distinctive profiles of all phase field solutions, 
which are almost constant-valued in each bulk region and vary smoothly but sharply in the thin 
diffuse interface. We refer to \cite{provatas2011phase, PDG99, du2008adaptive, elliott2010surface} for 
detailed discussions about adaptive simulations of 
phase field models from materials science and biology applications. It should be noted that
although coarse {\em a posteriori} error estimates were not rigorously derived in those cited works,  
they could be obtained by following the standard residual-based {\em a posteriori} 
estimate techniques \cite{EEJ95, ainsworth_oden00, nochetto_siebert_veeser09}. Moreover, we note that 
an alternative technique based elliptic reconstruction was also developed for deriving  
{\em a posteriori} error estimates for parabolic PDEs including the Allen-Cahn equation
\cite{lakkis_makridakis06, GLV11, LMP15}.

Like in the case of {\em a priori} error estimates for phase field models,  the coarse 
{\em a posteriori} error bounds is not useful anymore if one is interested in knowing 
the limiting behavior of the solution error and the error for numerical interfaces 
when $\varepsilon\to 0$. To this end, one must obtain fine {\em a posteriori} error estimates. 
As expected, the only known technique for deriving such error estimates is to utilize 
the PDE spectral estimate for the linearized operator of the underlying phase field PDE operator. 
Indeed, the desired fine {\em a posteriori} error estimates were obtained in
\cite{KNS04, Feng_Wu05} for the Allen-Cahn equation
and in \cite{Feng_Wu08} for the Cahn-Hilliard equation. We note that 
the arguments and the usages of the PDE spectral estimate are quite different 
in \cite{KNS04} and in \cite{Feng_Wu05, Feng_Wu08}. 
A topological argument was used in \cite{KNS04}, in which the PDE 
spectral estimate is embedded, to obtain a fine {\em a posteriori} error bound; while 
a sharp continuous dependence estimate for the PDE solution was used, in which the 
PDE spectral estimate was crucially utilized, to derive fine {\em a posteriori} error estimates 
in \cite{Feng_Wu05, Feng_Wu08} (also see \cite{cockburn03}). Because this continuous 
dependence argument is very simple and 
may be applicable to other evolution PDEs, we briefly explain it below. 

Let $V$ be a Hilbert space and ${\mathcal L}$ be an operator from
$D({\mathcal L})$ ($\subset V$), the domain of ${\mathcal L}$, to $V^*$,
the dual space of $V$. Consider the abstract evolution problem
\begin{alignat}{2}\label{eqn1}
\frac{\partial u}{\partial t} + {\mathcal L}(u)&= g &&\quad\mbox{in }\Omega_T:=\Omega\times (0,T), \\
u(0) &=u_0 &&\quad\mbox{in }\Omega. \label{eqn2}
\end{alignat}
Assume that the above PDE problem satisfies the continuous dependence estimate
in the sense that if $u^{(j)}$ is the (unique) solution of \eqref{eqn1}-\eqref{eqn2}
with respect to the data $(g^{(j)}, u_0^{(j)})$ for $j=1,2$, then there holds 
\begin{equation}\label{eqn3}
\|u^{(1)}-u^{(2)}\|_{L^p(0,T; V)} \leq \Phi\bigl(g^{(1)}-g^{(2)}\bigr)
+ \Psi\bigl(u^{(1)}_0-u^{(2)}_0 \bigr) 
\end{equation}
for some (monotone increasing) nonnegative functionals $\Phi(\cdot)$ and $\Psi(\cdot)$.
and $1\leq p\leq \infty$.  
Let $u^A$ be an approximation of $u$ with the initial value $u_0^A$, it is easy to show that
\cite{Feng_Wu05} 
	\begin{align}\label{eqn4}
	\|u-u^A\|_{L^p(0,T;V)} \leq \Phi\bigl(R(u^A) \bigr) + \Psi\bigl(u_0-u^A_0\bigr),
\end{align}
where
\[
R(u^A):= g-\frac{\partial u^A}{\partial t} - {\mathcal L}(u^A) .  
\]

Clearly, $R(u^A)$ denotes the residual of the approximation $u^A$. To utilize the above 
continuous dependence result to derive {\em a posteriori} error estimates for a spatial numerical discretization
method, one only needs to set $u^A = u_h$, the space-semi-discrete numerical solution and 
to get an upper bound for the residual $R(u_h)$ usually in terms of the element-wise local 
residuals and the jumps of the flux across the interior element edges/faces (cf. \cite{Feng_Wu05, Feng_Wu08}).
The biggest advantage of this approach for {\em a posteriori} error estimates is that it converts a 
numerical error estimate problem into a PDE continuous dependence (or stability) problem, which may be 
easier to cope with, especially, for deriving fine {\em a posteriori} error estimates for phase field models. 
In that case, the functionals $\Phi(\cdot)$ and $\Psi(\cdot)$ are expected to grow in $\frac{1}{\varepsilon}$
polynomially as shown in \cite{Feng_Wu05, Feng_Wu08}. 

Finally, we point out that the spectral estimate idea has also been successfully used as 
an {\em a posteriori} indicator 
for detecting singularities (such as topological changes) of the underlying sharp interface limit 
of a phase field model in \cite{antil2017spectral, Bartels15, Bartels_Muller11a, Bartels_Muller11b, 
	Bartels_Muller_Ortner11, Bartels_Muller10}. The rationale for this approach is that if no 
singularity occurs,  the 
discrete principal eigenvalue of the linearized operator, which is computable, should have a finite negative low bound. Hence, 
when such a low bound ceases to exist (i.e., it tends to $-\infty$ faster than $O(|\ln \varepsilon|)$ order 
as $\varepsilon\to 0$), then it flags a possible singularity of the underlying moving interface problem, 
and suggests that a fine mesh should be used to resolve the singularity.

\section{Applications and Extensions}\label{sec-7}

The phase field (diffuse interface) method has many applications ranging from  mathematical subjects like  
differential geometry, to image processing and geometric modeling,  and  to physical sciences like  astrophysics, 
cell biology, multiphase fluid mechanics, and (of course) materials science. We provide a sampler of some of the applications here.

\subsection{Materials Science Applications}\label{sec7.1}

Applications to problems in materials science are among the early and most widely recognized successes of 
phase field models. From the works of Lord Rayleigh, Gibbs and Van de Walls to the seminal contributions of John W. Cahn, 
ideas of using diffuse interface  and phase order parameters have drawn much interests in the materials science 
community in the mesoscopic modeling of materials structure.  Generically, a set of conserved field variables 
$c_1, c_2, ..., c_m$ and non-conserved field variables $\eta_1,\eta_2,...\eta_n$ are often used to describe 
the compositional/structural domains and the interfaces, and the total free energy of an inhomogeneous 
microstructure system is formulated as 
\begin{align}\label{pfenergy}
  E_{total} &=\int [\sum^m_{i=1}\alpha_i(\nabla
        c_i)^2+\sum^3_{i=1}\sum^3_{j=1}\sum^n_{k=1}\beta_{ij}\nabla_i\eta_k\nabla_j\eta_k \\
        &+ f(c_1, c_2, ..., c_m,\eta_1,\eta_2, ...\eta_n) ]
        dx+\int\int G(x-x', \vec{c}, \vec{\eta})dx dx' \;, \nonumber
\end{align}
where the gradient coefficient $\alpha_i$ and $\beta_{ij}$ can be used to reflect the interfacial energy anisotropy and the function $f$ corresponds to the local free energy density. The last integral in the above equation represents a nonlocal term that includes a general long-range interaction such as elastic interactions in solids. The time evolution is governed by either the non-conserved (Allen-Cahn type for $c_j$'s) or conserved dynamics (Cahn-Hilliard type for $\eta_j$'s).

Phase field modeling is arguably one of the most popular methods for modeling and simulation of microstructure evolution under different 
driving forces such as compositional gradients, temperature, stress/strain, and electric and magnetic fields \cite{biner2017}.  
Phase field modeling has now been used to study not only binary phases and single components but also materials with multicomponents multiphases. 
It also can account for mechanical, thermodynamic, electric and magnetic interactions. It can deal with  spinodal decomposition,  dendritic growth, 
Ostwald ripening and coalescence, adhesion and dewetting,  sintering, directional solidification, nucleation and coarsening, grain boundary 
motion, pattern formation in thin films, epitaxial growth, electromigration and  many other processes. 
For additional references and more comprehensive reviews, we refer to 
\cite{bgn1,bgn2,biner2017,bn1,bn2,bgn4,bgn6,bns,boettinger2002phase,chen2002phase,dai2016computational,provatas2011phase,SAL07,gns2004,singer2008phase,steinbach2009topical,wang2010phase}.

As an example, we mention one application on the study of nucleation based on the phase field approach \cite{zhang2016recent}. The search for index-1 saddle point of the phase field free energy \cite{zhang2007morphology,heo2010incorporating} can help  characterize the transition states that lead to critical nuclei. For solid state phase transformation, anisotropic elastic energy plays a critical role, which leads to formation of nonspherical nuclei. A particular form of the total energy in the phase field set up  is given by
$$
E_{total}= \int (| \nabla \eta|^2 + f(\eta)) dx + \frac{1}{2} \int C_{ijkl}
C_{ijkl}\varepsilon^{el}_{ij}\varepsilon^{el}_{kl} dx,
$$
where the elastic strain $\varepsilon^{el}$ is the difference between the total strain and stress-free strain since stress-free strain does not 
contribute to the total elastic energy. In  \cite{zhang2007morphology,zhang2008diffuse}, different critical nuclei have been computed as the 
driving forces change and the elastic effects vary. Numerical algorithms for index-1 saddle points are constructed by essentially reformulating 
the problem into a problem of optimization in an extended configuration space which includes not only the state variable (point on the energy 
landscape) but also the most unstable direction. In other words, one is effectively doing minimization along the stable directions (descending) 
while maximization along the (most) unstable direction (ascending) directions. This is the essential behind various saddle-point algorithms such as  
the Shrinking Dimer Methods \cite{zhang2012shrinking} and other methods reviewed in \cite{zhang2016recent}. The predicted phase field nuclei can 
also be utilized in dynamic simulations and extended to investigate complex nucleation phenomena \cite{heo2010incorporating}.

Let us also mention other variations of phase fiels models related to
materials science applications. For example, phase field modeling
of  polymeric materials has received much attention in theoretical
chemistry and materials science. Phase field models
of diblock copolymers, such as the model due to
 Ohta and Kawasaki \cite{ohta86}, have had a long history and is an active
area of research \cite{choksi2003derivation,ren2007many,wanner2017computer,ceniceros2004numerical,shi2019self,li-zhang_zhang13,cheng2017efficient}.

In addition, as noted previously, there are many other extended phase field models but with vector or tensor valued fields as phase field variables, 
These models can describe physical phenomena and geometric objects with higher co-dimension such as point defects  in two dimension and curves in 
three dimension. 

\subsection{Fluid and Solid Mechanics Applications}\label{sec7.3}

Applications of phase field models and diffuse interface approaches to mechanical problems is also a major area of research. 
 
Concerning fluids mechanics,  a review was given in \cite{anderson_mcfadden_wheeler98} on phase-field/diffuse-interface models
 of hydrodynamics and their application to a wide variety of interfacial phenomena involving fluid flows.  
 Phase field models have been developed in various cases of fluid flows with  a length scale commensurate 
 with the diffuse interfacial width. Examples include small-scale flows and multiphase flows. The latter 
 subject, typically involving breakup and coalescence  such as fluid jets and droplets as well as fluid 
 mixing and other interfacial phenomena, has become particularly important in microfluidics and micro-process
  engineering, see \cite{prosperetti2009computational,worner2012numerical} for reviews on the various 
  numerical approaches. Phase field modeling of multiphase flows has been developed by a number of groups 
  along with the development and numerical analysis, see, for instance  \cite{abels2012thermodynamically,lowengrub1998quasi,badalassi2003computation,bavnas2017numerical,diegel2017convergence,gns2000,feng2005energetic,grun2013convergent,liu2003phase} and the references cited therein.  Various phase field based simulations of drop formation in microfluidic channels have been carried out, see for example \cite{yue2004diffuse}.  Related optimal control problems and their simulations have been considered in \cite{hintermuller2012distributed,hintermuller2017fully}. Three-dimensional  phase field simulations of interfacial dynamics in Newtonian and viscoelastic fluids can be found \cite{zhou20103d}.  Algorithm development and numerical analysis of three phase flows can be found in \cite{yang2017numerical}.

Phase field modeling and simulation of multiphase flows with additional complications have also been made, studies include, for example, the wetting phenomena \cite{jacqmin2004onset}, two-phase flow in porous media \cite{cogswell2017simulation}.  A problem of considerable interests is the moving contact line problem for the interaction of fluid-fluid interface with a solid wall is a widely studied subject, see \cite{pomeau2002recent} for a review.  Ideas of diffuse interface and phase field models have also been considered \cite{jacqmin2000contact,pismen2001nonlocal}. A phase field model with generalized Navier boundary condition (GNBC)  is proposed in \cite{qian2003molecular}. 
Further computational studies can be found in \cite{HGW2011,BSSW2012,SBW2013,GW2014,ZW2016}.

Phase-field modeling has also received much attention in solid mechanics.   Representative works in this direction includes the modeling of 
brittle fracture \cite{borden2012phase,conti2016phase,miehe2010thermodynamically}. For example, in addtion to the fracture energy of the type 
given by the Cahn-Hilliard functional as proposed in \cite{bourdin2000numerical}, one can define the elastic energy density as
$$
\psi_e(\varepsilon, u) = [(1-k) u^2 + k] \psi^+_e(\varepsilon)
 k \psi^-_e(\varepsilon),
$$
where $u$ is the phase field variable,  $\psi^\pm_e$  are the strain energies computed from the positive and negative components of the strain 
tensor $\varepsilon^\pm$, respectively, defined through a spectral decomposition of strain \cite{miehe2010phase}, that is, 
$$\psi^+_e (\varepsilon)
=\frac{1}{2} \lambda[ ( {\rm Tr}  \varepsilon)^+]^2 
+ \mu {\rm Tr} [ ( \varepsilon^+)^2],
$$
$$\psi^-_e (\varepsilon)
=\frac{1}{2} \lambda[ {\rm Tr}  \varepsilon ( {\rm Tr}  \varepsilon )^+]^2 
+ \mu {\rm Tr} [ ( \varepsilon - \varepsilon^+ )^2].
$$
The phase-field variable is only applied  to the tensile part of the elastic energy density, so that crack propagation under compression
 can be prevented \cite{miehe2010thermodynamically,miehe2010phase}. One may further incorporate other features to model ductile fracture \cite{ambati2016phase,miehe2016phase} and fluid-filled fracture \cite{wheeler2014augmented},   and utilize various discretization techniques, 
 see recent discussions in \cite{heister2015primal,zhang2018study}.  

\subsection{Image and Data Processing Applications}

Ideas of phase field modeling are used not only  for physical sciences  but also other domains. Imaging science is a good example where 
variational methods have been widely used for imaging and data analysis \cite{chan2005image,aubert2006mathematical}. Phase field models 
and diffuse interface descriptions of geometric features form one major  type of variational formulations  for various image processing tasks. 
Taking, for example, the task of image restoration and image deblurring, a diffuse interface relaxation to Mumford-Shah can be formulated 
by \cite{brett2014phase}
\begin{align}\label{phasefield-deblurring} 
\mathcal{J}(u):= \int_\Omega \Bigl( \frac12 |\nabla u|^2 + \frac{1}{\varepsilon^2} 
 ( u^2-1 )^2 \Bigr) \, dx + \|Su - f\|^2,
\end{align}
where $S$ can be a deblurring operator. 
Other examples include image classification and restoration \cite{samson2000variational},  edge detection \cite{AT1990}, motion estimation \cite{preusser2007phase},  image segmentation \cite{BOS2013,jung2007multiphase,li2011multiphase},
curve smoothing \cite{zhu2010variational}, active contour \cite{rochery2005phase}, clustering of vector fields \cite{garcke2001phase}, and so on.

For image segmentation,
a phase field relaxation to Mumford-Shah 
can be formulated 
as  \cite{march1997variational} 
\begin{align}\label{phasefield-segment} 
\mathcal{J}(v,u):=  E(u) + \lambda 
\int_\Omega u^2 |\nabla v|^2 dx + 
\|v - f\|^2\,,
\end{align}
where $f$ is the image under consideration, $v$ represents a piecewise smooth approximation to $f$, and  $E(u)$ is a diffuse interfacial 
energy of the phase field variable $u$. 

While a standard Cahn-Hilliard energy is the most popular choice, a more general form has also used as well:
\begin{align}\label{phasefield-segment-1} 
E(u): = \alpha
\int_\Omega \Bigl( \frac12 |\nabla u|^2 + \frac{1}{\varepsilon^2} 
 ( u^2-1 )^2 \Bigr) \, dx + 
\beta
\int_\Omega \left|\Delta u+ \frac{1}{\varepsilon^2} 
 (u- u^3) \right|^2  \, dx.
\end{align}
The first term in the above is the standard phase field form of the surface tension, the second term corresponds to a  diffuse interface 
relaxation of the Euler-elastica (integral of the square of curvature along a planar curve) as discussed in \cite{de1991some}.

Discussions on issues related to the application of similar phase field energy to image inpainting can be found in \cite{esedoglu2002digital}.
 More discussions on the three dimensional analog of the phase field Euler-elastica, 
a special form of the Helfrich bending energy,  are given later in this chapter for phase field modeling of vesicle  membrane.

%graph
Phase field approaches have also been used in data analysis and discrete graph modeling, see for example 
\cite{bertozzi2012diffuse}, \cite{van2014mean}, and \cite{li2011multiphase}.
    
\subsection{Biology Applications}\label{sec7.2}
In more recent years, phase field modeling and simulations have also become very successful tools to study 
problems arising in  life sciences, including problems in biophysics, bioengineering, ecology, biology, and so on. 
Some examples include the development of phase field models and computational methods for the deformation 
and dynamics of bio-membranes and cell vesicles, cellular activities such as  motility and cell division, fluid-structure
 interactions in blood flows and other biological fluids, and tumor adhesion and growth.
These problems span across scales from cellular and molecular levels to macroscopic systems.

On the cellular level, activities such as cell crawling, cell motion on patterned substrate and other migration processes are also important biological questions. They often involve deformable geometry and associated mechanical interactions. These questions  have been studied using various phase field approaches in recent years, see for example \cite{camley2013periodic,lober2014modeling,marth2014signaling,moure2017phase,najem2016phase,shao2012coupling,ziebert2016computational}. Phase field modeling has also found many applications in biological molecular processes such as conformational change, molecular recognition, and molecular assembly. For example, \cite{sun2015self} presented phase field models for the implicit solvation of charged molecules with a coupling to Poisson--Boltzmann electrostatics. Various refined phase field models,numerical methods on both spatial discretization and time integration have also been studied \cite{LZ13,DLL18,Zhao2018}.

Bilayer vesicles are models systems for cell membranes. Their deformation has often been modeled by the Helfrich bending  elasticity model for fluid membranes with contributions from the surface integrals of the mean curvature squared and the Gaussian curvature, thus leading to natural connections with the Willmore energy and Willmore flow. The diffuse interface formulation has been introduced in \cite{du2004phase,du2006simulating}. Other variants have also been studied \cite{aland2014diffuse,biben2005phase,bretin2015phase,bgn7,campelo2007shape,du2009energetic,jamet2007towards,esedoglu2014colliding}, together withs analytical studies \cite{BMO15,Mu13,bm2010} and 
numerical approximations \cite{CL11,chen2015decoupled,du2006analysis,du2008adaptive,yang2017efficient} for the deformation and dynamics of vesicles. 
The effective phase field modeling of Gaussian curvature energy also leads to an interesting development of the diffuse interface Euler-Poincar\'e characteristics \cite{du2005retrieving,bm2010} that can be used to detect topological change of the implicitly defined interface within the phase 
field framework beyond the biophysical applications. Phase field formulation has also been developed for vesicle-substrate and vesicle-vesicle interactions \cite{zhang2009phase,gu2016two}. For multicomponent membranes, \cite{elliott2010surface,LRV09} studied phase field models defined on 
an evolving surface. Using two phase field variables, \cite{wang2008modelling} developed phase field models for two-component membranes and obtained
 interesting patterns mimicking to experimental observations. Given the complexities involved in these systems, effective numerical algorithms such 
 as those utilizing high order spectral methods and adaptive discretizations \cite{du2004phase,du2008adaptive}  and energy law preserving schemes 
 \cite{hua2011energy} can be very useful for large scale and long time simulations. 

Another example of broad applications of phase field modeling and simulations in biomedical field is on the  cancer and tumor growth models. 
This is a subject  that have  been  studied by a number of authors\cite{garcke2018multiphase,vilanova2017mathematical,lima2014hybrid},  see \cite{lowengrub2009nonlinear} for a review on the models and numerical methods.

\subsection{Other Variants of Phase Field Models}\label{sec7-3}

There are many other variants of phase field models. We provide a couple of examples here that focus on the modeling of the nonlocal and stochastic aspects of the underlying  physical processes.

\subsubsection{Nonlocal and Factional Order Phase Field Models}

Discussion on nonlocal interactions in the form of integral operators may be tracked back to the work of 
Van der Waals \cite{van1894thermodynamische}, see discussions made in \cite{pismen2001nonlocal}, The usual differential 
equation form of the local phase field energy can be derived from the nonlocal version via the so-called Landau 
expansion \cite{landau2013course}, assuming a smooth and spatially slowly varying field.  A number of studies on nonlocal 
Allen-Cahn and nonlocal Cahn-Hilliard can be found in \cite{bates2006nonlocal,bates2006some,benesova2014implicit,choksi2009phase,choksi20112d,djlq19,fife2003some,jeong2015microphase}. More studies on nonlocal modeling, analysis and computation can be found in \cite{du19}. A special class of nonlocal models are fractional phase field models where fractional derivatives are used to replace the 
integer order derivatives in the conventional local phase field models. Such models have attracted the  attention of community 
in recent years \cite{akagi2016fractional,gui2015traveling,milovanov2005fractional,tarasov2005fractional,valdinoci2013fractional}.
Algorithmic development and numerical analysis concerning these nonlocal models (including the fractional ones in either space 
or time or both ) can be found in \cite{ainsworth2017analysis,bates2009numerical,du2016asymptotically,du2018stabilized,guan2014second,hou2017numerical,song2016fractional,zhai2016fast}. Related algorithmic studies with respect to different applications were presented in  \cite{antil2017spectral}.

For spatially nonlocal phase field models, due to the spread of nonlocal interactions, the phase field variables may no longer be 
as smooth as their local counterpart. They may also lead to narrower interfacial region and permit singularities across the interface 
or at the defects \cite{du2016asymptotically,gui2015traveling,song2016fractional}.

\subsubsection{Stochastic Phase Field Models}

Uncertainty may arise from various sources such as thermal fluctuation, impurities of the materials and the 
intrinsic instabilities of the deterministic evolutions. Therefore, the evolution of interfaces under influence
of noise is of great importance in applications such as lattice models, the scaling limit of lattice models and 
derivations of continuum equations, it is necessary and interesting 
to consider stochastic effects, and to study the impact of noise on phase field modeling and on their solutions, especially 
on their long time behaviors. 
This then leads to considering the stochastic phase field models. However, how to incorporate noises correctly into 
those models is a nontrivial matter, which turns out is both a science and an art. 

A few approaches have been known in the literature. The first one, which
is the simplest, is to add (small) noise terms to the existing deterministic models and to study and simulate those 
stochastically perturbed models. Some recent works on modeling and PDE analyses in this direction can be found in
\cite{E_liu02, Weber10a, Weber10b, hrw12, OWW14,ABK2018,afk18,DZ07,DG11,DPD96,RYZ18} and the references therein.
Finite element approximations have been studied in \cite{kovacs2011finite, KLM14, FKLL18} for the Cahn-Hilliard-Cook equation 
(which is a stochastic Cahn-Hilliard equation with additive noise) and in \cite{KLL15, FLZ18} for stochastic Allen-Cahn models. Relevant time discretizations 
have been discussed in \cite{printems2001discretization}. Parallel algorithms and numerical simulations of the Cahn-Hilliard-Cook 
equation have also been reported in \cite{zheng2015parallel}. Stochastic Cahn-Hilliard dynamics have also been used to study 
nucleation in microstructure evolution in \cite{heo2010incorporating,li2012numerical,li-zhang_zhang13}. It should be noted that the added 
noises may not always have physical meaning and those stochastically perturbed models may not associate with sharp interface models. 

An alternative  approach  is to directly consider stochastic sharp interface problems such as stochastic mean 
curvature flows and to derive (formally) corresponding phase field models, subsequently, to study and simulate 
the resulting stochastic models. Recent PDE studies in this direction can be found in \cite{weber_roger13, dln01, souganidis_yip04,Yip02,Yip98}, those stochastic PDEs involve gradient-type multiplicative noise and thus have stronger nonlinearity. 
Finite element approximations have been carried out in \cite{FLZ17} for a stochastic Allen-Cahn equation with 
gradient-type multiplicative noise and in \cite{FLZ19} for a related stochastic Cahn-Hilliard equation with 
gradient-type multiplicative noise. Not surprisingly, their sharp interface models are expected to be a stochastic 
mean curvature flow \cite{weber_roger13} and a stochastic Hele-Shaw model \cite{FLZ19, abk18}, respectively. However, 
their rigorous convergence proofs are still missing although some partial results were reported in \cite{weber_roger13, abk18}
and positive numerical results were given in \cite{FLP14,FLZ17,FLZ19}. 

An additional approach to model stochasticity is
to study  stochastic variants of coupled phase field 
 equations with other physical models. 
Stochasticity can enter through other physical laws such as equations
for fluid flows in multiphase flows and flow-structure interactions
\cite{chaudhri2014modeling,du2011analysis}.

\section{Conclusion}\label{sec-9} 
In this chapter, we presented a holistic review about some basic elements of phase field modeling, 
particularly related to curvature-driven geometric interfacial motion.  We discussed the relevant mathematical 
theory, numerical approximations and selected applications, as well as the relationship between 
the phase field methodology and other methodologies (such as the level set methodology)
for geometric moving interface problems. Instead of presenting much involved technical details, 
we focused on discussing the main elements and ideas of the phase field modeling, analysis and their numerical 
approximations, with extensive references provided  on each of these aspects. 
They should be helpful for the interested reader to do further reading for specific details. 

Over last several decades, the phase field method has been developed into a powerful and versatile 
general methodology for interface problems arising from various scientific and engineering 
applications including biology, differential geometry, fluid and solid mechanics, image processing
and materials science. It has garnered tremendous attention and popularity among researchers and 
practitioners, the trend will likely continue and make even broader impact in more and more fields.
Yet, there are still many challenging issues in the phase field methodology, from modeling, analysis, approximation
and application, that should
be further investigated. On the modeling side, addressing different geometric features and connecting to microscopic
physics are important future research topics. Effective recovery of geometric and topological features 
and statistical information from the phase field approach are of both theoretical and practical interests.
On the computational side, there are much need and urgent demand for developing 
effective and robust adaptive algorithms, particularly those involving anisotropic spatial adaptivity 
and locally adapted time steps. Developing efficient and fast linear and nonlinear solvers,  
preconditioning techniques, high order and stable time stepping schemes, and scalable algorithms 
remains a focus of future research. Effective methods for computing saddle 
points or transition states, rather than equilibrium, are also very limited currently and are in high demand.
Methods for characterizing the captured interfacial geometry and quantifying statistical features based on the
 phase field models are also helpful for many practical applications.
On the analysis front, a lot of open questions remain to be answered. For the PDE analysis,  convergence 
studies of many phase field models to their respective sharp interface limits are still missing.  The lack of 
analysis techniques and machineries along with singularities of solutions of the sharp interface limits
seem to be one of the primary hurdles to overcome. For the numerical analysis, analyzing the approximation 
errors and their dependences on the diffuse interface parameter $\varepsilon$ remains an interesting 
and challenging issue.
Finally, while extreme scale simulations have been achieved for some model phase field equations, 
simulating more realistic and more complex phase field systems still require much future effort.

%\newpage
\nocite{*}
 
%%%%%%%%%%%%%%%%%%%%%%%%%%
\bibliographystyle{abbrv}
%\bibliographystyle{plain}
%\section*{References}

%\bibliography{./references}
%%%%\bibliography{./newref}

\end{document}